\documentclass[a4paper,11pt]{article}
\usepackage{latexsym,amssymb,enumerate,amsmath,epsfig,amsthm}
\usepackage[margin=1in]{geometry}
\usepackage{setspace,color}
\usepackage{tikz}
\usepackage{floatrow}
\usepackage{subfig}
\usepackage[ruled]{algorithm2e}

\newcommand{\p}{\mathbf{p}}
\newcommand{\e}{\mathbf{e}}
\newcommand{\X}{\mathbf{X}}

\newcommand{\Surface}{\Sigma}

\newcommand{\reminder}[1]{{#1}}

\begin{document}

\title{Modified Virtual Grid Difference for Discretizing the Laplace-Beltrami Operator on Point Clouds}
\author{
Meng Wang\thanks{
              Department of Mathematics,
            Hong Kong University of Science and Technology,
            Clear Water Bay, Hong Kong ({mwangae@ust.hk},{masyleung@ust.hk}).}
\and
Shingyu Leung\footnotemark[2]
\and
Hongkai Zhao\thanks{
            Department of Mathematics,
            University of California at Irvine,
            Irvine, CA92697-3875, USA ({zhao@math.uci.edu}).}
}
\maketitle

\begin{abstract}
We propose a new and simple discretization, named the Modified Virtual Grid Difference (MVGD), for numerical approximation of the Laplace-Beltrami (LB) operator on manifolds sampled by point clouds. The key observation is that both the manifold and a function defined on it can both be parametrized in a local Cartesian coordinate system and approximated using least squares. Based on the above observation, we first introduce a local virtual grid with a scale adapted to the sampling density centered at each point. Then we propose a modified finite difference scheme on the virtual grid to discretize the LB operator. Instead of using the local least squares values on all virtual grid points like the typical finite difference method, we use the function value explicitly at the grid located at the center (coincided with the data point). The new discretization provides more diagonal dominance to the resulting linear system and improves its conditioning. We show that the linear system can be robustly, efficiently and accurately solved by existing fast solver such as the Algebraic Multigrid (AMG) method. We will present numerical tests and comparison with other exiting methods to demonstrate the effectiveness and the performance of the proposed approach.
\end{abstract}

\section{Introduction}
Solving partial differential equations (PDEs) on surfaces or manifolds in general has many important applications arising from problems in science and engineering. Examples include diffusion process simulation on interfaces in multiphase problems in fluids and materials, three dimensional geometric modeling and shape analysis, manifold learning in high dimensional data analysis, to name just a few. One numerical representation of the surface is given by a point cloud which is the simplest, the most natural and ubiquitous way of sampling and representing surfaces and manifolds of free form in three dimension and higher. Points cloud data are routinely obtained by modern sensing technology and are extensively used in 3D modeling and shape analysis. Feature vectors used in data analysis can also be viewed as point clouds typically embedded in a high dimensional Euclidean space.

Various intrinsic differential operators have been successfully used to \textit{connect the dots} to extract geometric quantities and local/global structures from these unstructured point clouds.  Hence, it is desirable to develop efficient and robust numerical methods for solving PDEs directly on point clouds without a global mesh or parameterization. Although point clouds provide a simple and flexible way for geometric representations, it leads to difficulties when computing integrals and solving differential equations. In this work, we propose a simple numerical method for discretizing the Laplace-Beltrami (LB) operator on manifolds sampled by point clouds.  The LB operator is one of the most important differential operators defined on manifolds. It models diffusion process on surfaces in physics. The eigen-system of LB operator also provides an intrinsic orthogonal basis for square integrable functions and it contains intrinsic geometrical information of the underlying manifold. The LB operator and related tools such as the LB eigen-map, the heat kernel and diffusion maps are used extensively in shape analysis and manifold learning.

Given a manifold sampled by a point cloud, there are a few ways to approximate the LB operator.
One way is first to construct a global mesh, e.g., triangulation, on the point cloud. Then the LB operator can be readily approximated using any standard finite element method based on variational formulation on the triangulation \cite{dzi88,tau00,mdsb03,xu04}. This approach naturally results in a symmetric positive definite linear system, which is the main advantage of the approach.  On the other hand, obtaining a global triangulation itself can already be a challenging task. Moreover, the resulting triangulation may not be good enough for a stable finite element method if the sampling of the point cloud is non-uniform. Instead of building a global mesh, a local mesh method has been recently proposed in \cite{lailiazha13}. However, the discretized system is not symmetric in general.

There are also a few methods for approximating the LB operator directly on point clouds. Instead of building a mesh on point clouds, these methods only utilize neighboring points relation and are particularly useful when a good quality global mesh is intractable. One type in this class is the kernel based methods \cite{belniy05,coilaf06,belsunwan09,lishisun15,sinwu15}. The idea is to transform the differential equation into an integral equation. If the kernel function is the Green's function for the LB operator on the underlying manifold, an exact integral equation can be derived. However, since the Green's function for the LB operator is impossible to obtain in practice, one usually derive an approximate integral equation using a radial symmetric kernel function in the ambient Euclidean space with compact support or \reminder{exponential decay}. Even though kernel based method can be generalized to more general settings such as graph Laplacian and diffusion maps, the method is still mainly restricted to the diffusion operator. Kernel based methods typically lead to a discrete linear system with a M-matrix, which gives the discrete system a maximum principle analogous to that of the continuous problem. The approximation error is determined by two competing scales. One is the spatial scale of the kernel at which geodesic distance of the underlying manifold is approximated by the Euclidean distance in the ambient space. The other is the finer scale at which the data is sampled and such that integration at the kernel scale can be approximated accurately enough. Although one can show that the solution to the integral equation converges to the corresponding solution to the LB equation under certain sampling conditions of the point cloud, the order of convergence is low \cite{lishisun15,llwz12}.

Recently, motivated by our earlier work \cite{leuzha0801} where the grid based particle method (GBPM) was developed for moving interface problems, \cite{llwz12,liazha13} has introduced another framework for solving PDEs directly on point clouds. The GBPM represent and track a moving interface by meshless Lagrangian particles based on an underlying Eulerian mesh. This results in a special type of point clouds representation of the interface with a quasi-uniform sampling rate. We refer interested readers to \cite{leuzha0801,leuzha0802,leuzha10,leulowzha11,liuleu13,honleuzha14} for a complete discussion on the representation. In a follow-up work \cite{leulowzha11}, the method has been further developed to solve several kinds of PDEs on surfaces arising from physical and geometric flow based on the GBPM. The key observation is that both the manifold and the function on it can be parametrized in a local Cartesian coordinate system at each point of the point cloud. Hence, differential operators can be approximated at each point by applying a least squares approximation to both the manifold (and its metric) and the function in the local coordinate system through neighboring points. The least squares approximation in the approach introduce a more robust and flexible discretization compared to the exact interpolation on point clouds. The method applies to general differential operators on point clouds sampling manifolds with arbitrary dimensions and co-dimensions. Local least squares approximation allows also high order accuracy. Moreover, the computational complexity depends mainly on the intrinsic dimension of the manifold rather than the dimension of the embedded space. However, the computational cost for least squares approximation grows quickly with the dimension of the manifold. Although local approximation accuracy is relatively easy to achieve, a more important and more challenging issue is how to construct a discretization on the whole point cloud so that the resulting linear system can be solved efficiently and stably. For example, when applied to the LB operator, one can obtain a M-matrix system by using a constrained quadratic programming optimization technique to enforce both the consistency and diagonal dominance after the discretization \cite{liazha13}.

In this work, we propose a new and simple discretization of the LB operator on point clouds, called the Modified Virtual Grid Difference (MVGD). There are two key ideas in our new approach. Based on the same observation that both the manifold and the function on it can be parametrized in a local Cartesian coordinate system centered at each point of the point cloud, we first introduce a virtual grid with a scale adapted to the sampling density in the local coordinate system and we apply finite difference scheme on the virtual grid to discretize the LB operator. Secondly, functions values on the virtual grid used for the finite difference scheme are interpolated from the local least squares approximation of the function {\bf except} at the center grid (coincided with the data point), where the function value at the data point is explicitly used. As will be shown later, this MVGD discretization leads to a more diagonal dominant and a better conditioned linear system which can be efficiently and robustly solved by many existing fast solver for elliptic PDEs.

Using a simple 1D example with some explicit formula, we will show that the MVGD improves both the numerical accuracy and computational stability as compared to the original least squares approach.  The numerical implementation is also simpler compared to the discretization using the constrained quadratic optimization approach proposed in \cite{liazha13}. As demonstrated by numerical tests, the MVGD can handle point clouds with non-uniform sampling reasonably well too.

\reminder{It is worth mentioning that although in this paper we are only concentrated only on solving LB equation on point clouds, there have actually been some interesting recent research works on surface information approximation, parametrization, reconstruction and rendering \cite{frenie80,GSC95,kol98,krsek1998algorithms,meek2000surface,lazmon02,alexa2003computing,medvelfig03,pauly2003shape,fleishman2005robust,magsolriv07,liu2008local}. They provide alternative tools to provide a local approximation of the underlying manifold which is also important for the subsequent discretization of any differential operator on the underlying manifold.}

The paper is organized as follows. In Section \ref{Sec:Background}, we first briefly summarize the approach discussed in \cite{liazha13}, and will discuss some issues associated to the discretization of LB operator. The proposed method will be given in Section \ref{Sec:ModifyApp}. Explicit discretization for the LB operator for one and two dimensional manifolds will be explicitly constructed. A simple numerical examples in one dimension will be provided a few insights. In Section \ref{Sec:Examples}, we first test our new discretization by solving the LB equation on point clouds. We then compute the eigenvalues and the eigenvectors of the LB operator and make comparisons with other methods.


\section{Background: Approximating the LB operator on point clouds}
\label{Sec:Background}

For simplicity, we consider a two-dimensional surface $\Sigma$ in $\mathbb{R}^3$ parametrized by $(s_1,s_2)$, the LB  operator acting on a function $U:\Sigma \rightarrow \mathbb{R}$ is defined by
\begin{equation}
\Delta_{\Sigma} U= \sum_{{\alpha,\beta}=1}^2 \frac{1}{\sqrt{g}}
\frac{\partial}{\partial s_{\alpha}} \left( \sqrt{g} g^{\alpha \beta} \frac{\partial
U}{\partial s_{\beta}} \right) \label{Eqn:SurLapEqn}
\end{equation}
where the metric $[g_{\alpha \beta}]$ is given by
$$
g_{\alpha \beta}=\frac{\partial \X}{\partial s_{\alpha}} \frac{\partial \X}{\partial s_{\beta}} \, ,
$$
with the surface $\Sigma$ given by $\X(s_1,s_2)$. The function $g=g_{11}g_{22}-g_{12}g_{21}$ is the Jacobian of the metric and $[g^{\alpha \beta}]$ is the inverse of the matrix $[g_{\alpha \beta}]$. \reminder{Note further that the Laplace-Beltrami operator is geometric intrinsic and is therefore independent of the parametrization. So our later computations based on local parametrization is still valid and geometric intrinsic.}

Numerically, given a point cloud $S=\left\{ \mathbf{p}_i \mid i=1,...,N \right\}$ sampling $\Sigma$, we want to approximate the above LB operator at any point $\mathbf{p}_i\in S$. A general framework has been proposed in \cite{liazha13}. The key idea is that both $\Surface$ and $U$ are functions that can be parametrized in a local Cartesian coordinate system centered at $\p_i$, and can be approximated to any desired order, e.g., by polynomials using least squares on neighbors of $\p_i$. Then the LB operator \eqref{Eqn:SurLapEqn} at $\p_i$  can be defined in the local coordinate system and easily approximated as well.  A key point is that the resulting numerical approximation or discretization of the LB operator of a  function at a point $\p_i$ is written as a linear combination of the function values at neighboring points of $\p_i$. The coefficients of the linear combination depends on the relative locations of the neighboring points and the choice of least squares approximation. Although local approximation accuracy is relatively easy to achieve, the central issue is essentially how to construct a discretization so that the resulting linear system can be efficiently inverted in a stable manner. For example, it is desirable for the discrete system to preserve analogous properties of the underlying continuous differential operator. Without a global mesh, symmetry of the linear system is not possible in general. Even so, it is still possible and desirable to obtain a M-matrix for the linear system which gives discrete maximum principle and hence stability for the numerical solution. Here we give a brief summary of the idea. For a detailed description, we refer interested reader to \cite{liazha13} and some references thereafter.

For a given $\mathbf{p}_i$, we first collect its $K$ neighboring points, e.g., $K$ nearest neighbors (KNN), and denote them by $N(i)=\{\p_j: \p_j$ belongs to the nearest $K$ neighboring points of $\p_i$, for $j=1,2,\cdots,K$\}. Then the manifold is locally approximated using these points in an appropriate Cartesian coordinate system centered at $\p_i$. For example, one can use standard Principal Component Analysis (PCA) to obtain the local coordinate system. Let $C_i$ be the covariance matrix at point $\mathbf{p}_i$ through its KNN,
$$
C_i = \sum_{\mathbf{p}_k\in N(i)}(\mathbf{p}_k-\mathbf{c}_i)^T(\mathbf{p}_k-\mathbf{c}_i)
$$
where $\mathbf{c}_i$ is the local barycenter of $N(i)$ given by
$$
\mathbf{c}_i = \frac{1}{K}\sum_{\mathbf{p}_k\in N(i)}\mathbf{p}_k \, .
$$
The eigenvectors of the covariance matrix $C_i$ provides a Cartesian coordinate system centered at $\mathbf{p}_i$. To simplify the later notations, we denote the corresponding coordinates of KNN by \{$(x_j,y_j,z_j), \p_j \in N(i)$\}. Here the $z$-axis is in parallel to the eigenvector that corresponds to the smallest eigenvalue and is an approximation of the surface normal at $\p_i$.

\reminder{As long as the surface is smooth, in a sufficient small neighborhood of $\mathbf{p}_i$ resolved by local sampling density, we can construct a local Cartesian coordinate system in which the surface $\Sigma$ is a graph.} Then it can be locally approximated by a polynomial, e.g., of second order, denoted by $z^{(i)}$
\begin{equation}
z^{(i)}(x,y)= \sum_{\alpha=0}^2 \sum_{0 \le \alpha+\beta \le 2} a^{(i)}_{\alpha,\beta}
x^{\alpha} y^{\beta} \, ,
\label{eq:poly}
\end{equation}
where the metric $g_{\alpha \beta}$ in this coordinates system can be approximated by
\begin{equation}
[g_{\alpha \beta}]=\left(
\begin{array}{cc}
1+\left(\frac{\partial z^{(i)}}{\partial x}\right)^2 &
\frac{\partial z^{i)}}{\partial x}\frac{\partial
z^{(i)}}{\partial y} \\
\frac{\partial z^{(i)}}{\partial x}\frac{\partial z^{(i)}}{\partial y} & 1+\left(\frac{\partial z^{(i)}}{\partial y}\right)^2
\end{array}
\right) \, ,
\end{equation}
with
\begin{eqnarray}
\frac{\partial z^{(i)}}{\partial x} &=& a^{(i)}_{1,0}+2a^{(i)}_{2,0} x +
a^{(i)}_{1,1} y \, , \nonumber\\
\frac{\partial z^{(i)}}{\partial y} &=& a^{(i)}_{0,1}+2a^{(i)}_{0,2} y +
a^{(i)}_{1,1} x \, .
\label{eq:metric}
\end{eqnarray}

\begin{figure}[!htb]
\begin{center}
\includegraphics[width=8cm]{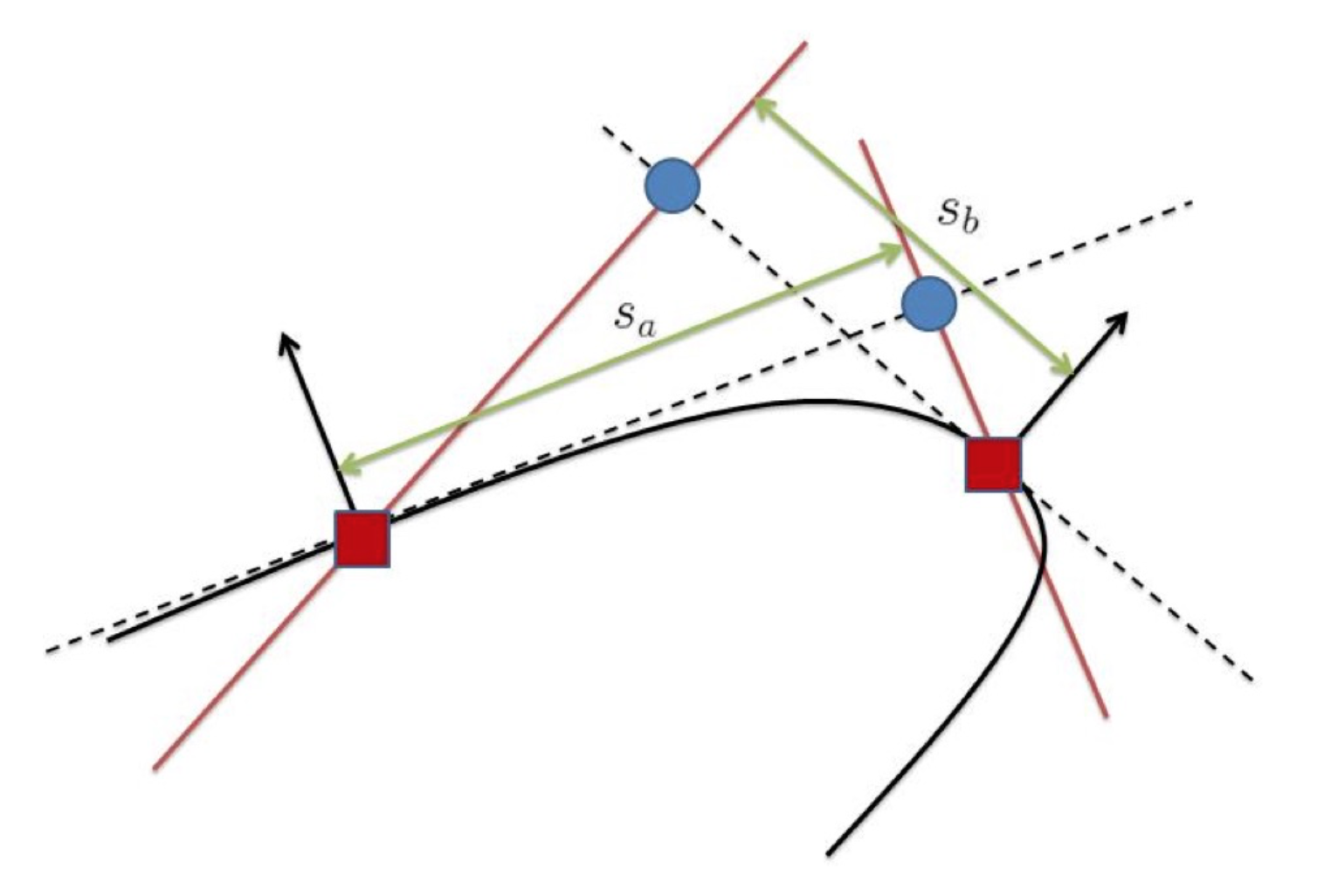}
\end{center}
\caption{Discretization of the surface Laplacian gives a nonsymmetric matrix.} \label{Fig:UnsymMat}
\end{figure}

Numerically, the local approximation \eqref{eq:poly} can be computed using (weighted) moving least squares (MLS). Given a function $U$ defined on $\Sigma$ also sampled on the point cloud, we again use MLS to approximate it in the same local coordinate system at $\p_i$, e.g., by a quadratic polynomial,  denoted by $U^{(i)}$,
\begin{equation}
U^{(i)}(x,y)= \sum_{\alpha=0}^2 \sum_{0 \le \alpha+\beta \le 2} b^{(i)}_{\alpha,\beta} x^{\alpha} y^{\beta}
\, .
\label{eq:T}
\end{equation}
Note that such least squares approximation is valid only locally near $\p_i$. \reminder{For general situations when the number of data points is larger than the degree of freedom or when the underlying function $U$ is not a polynomial of the same or lower degree, the value of least square approximation at $\p_i$ given by $U^{(i)}(0,0)$ is generally different from $U(\p_i)$. Once we have obtained such least squares approximation, we can then approximate the derivatives of $U$ using the corresponding derivatives of the polynomial $U^{(i)}(x,y)$}. By expanding (\ref{Eqn:SurLapEqn}), we obtain
\begin{equation}
\Delta_{\Surface} U^{(i)}(x,y) \simeq
A_1^{(i)} \, \frac{\partial U^{(i)}}{\partial x}+A_2^{(i)} \, \frac{\partial U^{(i)}}{\partial y} + A_3^{(i)} \, \frac{\partial^2 U^{(i)}}{\partial x^2}+A_4^{(i)} \, \frac{\partial^2 U^{(i)}}{\partial x \partial y}+A_5^{(i)} \, \frac{\partial^2 U^{(i)}}{\partial y^2} \, ,
\label{eqn3}
\end{equation}
where coefficients $A^{(i)}_l, l=1, \ldots, 5$ depend on the local metric, which can be computed from local reconstruction of the surface \eqref{eq:poly} to \eqref{eq:metric}. From the local approximation \eqref{eq:T} for $U^{(i)}$, we have:
$$
\frac{\partial U^{(i)}}{\partial x}=b^{(i)}_{1,0},~~~\frac{\partial U^{(i)}}{\partial y}=b^{(i)}_{0,1},~~~\frac{\partial^2 U^{(i)}}{\partial x^2}=2b^{(i)}_{2,0},~~~\frac{\partial^2 U^{(i)}}{\partial x\partial y}=b^{(i)}_{1,1},~~~\frac{\partial^2 U^{(i)}}{\partial y^2}=2b^{(i)}_{0,2} \, .
$$
As $U^{(i)}(x,y)$ is the least squares approximation,
we have the system of linear equation $\mathbf{A} \mathbf{b} = \mathbf{U}$ where
$$
\mathbf{A}= \left( \begin{array}{cccccc}
K+1 & \sum x_j & \sum y_j & \sum x^2_{j} & \sum x_{j}y_{j} & \sum y^2_{j} \\
\sum x_j & \sum x^2_j & \sum x_{j}y_{j} & \sum x^3_{j} & \sum x^2_{j}y_{j} & \sum x_jy^2_{j} \\
\sum y_j & \sum x_{j}y_{j} & \sum y^2_j & \sum x^2_{j}y_{j} & \sum x_{j}y^2_{j} & \sum y^3_{j} \\
\sum x^2_j & \sum x^3_j & \sum x^2_{j}y_{j} & \sum x^4_{j} & \sum x^3_{j}y_{j} & \sum x^2_{j}y^2_{j} \\
\sum x_{j}y_{j} & \sum x^2_{j}y_{j} & \sum x_{j}y^2_{j} & \sum x^3_{j}y_{j} & \sum x^2_{j}y^2_{j} & \sum x_{j}y^2_{j} \\
\sum y^2_j & \sum x_{j}y^2_{j} & \sum y^3_j & \sum x^2_{j}y^2_{j} & \sum x_{j}y^3_{j} & \sum y^4_{j}
\end{array} \right) \, ,
$$
$$
\mathbf{U}= \left(
\sum U_j,
\sum x_jU_j,
\sum y_jU_j,
\sum x^2_jU_j,
\sum x_jy_jU_j,
\sum y^2_jU_j\right)^T \, ,
$$
$\mathbf{b}=\left(
b^{(i)}_{0,0},
b^{(i)}_{1,0},
b^{(i)}_{0,1},
b^{(i)}_{2,0},
b^{(i)}_{1,1},
b^{(i)}_{0,2}\right)^T$ with the sum is taken over $\p_i$ and its KNN, i.e. $\{\p_j\in N(i)\}$. So, $b^{(i)}_{\alpha,\beta}$ for $\alpha,\beta=0,1,2$ are all linear combinations of values of $U$ at $\p_i$ and $\p_j\in N(i)$. Thus, at every point on the point cloud, the LB operator can be discretized and expressed as a linear combination of its value and its neighboring values, where the coefficients depend on the  locations of neighboring points. By assembling the discretization at all points of the point cloud together, we can obtain a linear system with the function value at each point of the point cloud as unknowns. Each row of the matrix corresponds to the discretization of LB operator at a point and we denote the matrix by $A$. First of all, $A$ is not symmetric in general because the least squares approximation at each point is computed in a local coordinate system through its nearest neighbors. For a general point cloud, the relative positions of nearest neighbors for two nearby points are different. This fact is illustrated by Fig \ref{Fig:UnsymMat}. Consider the two points (red square) on a curve (black curve). Local coordinate systems are plotted using arrow (along the normal) and dashed (along the tangent) lines, respectively. The $x$ coordinate of each point translated to the local coordinate system of the other point is denoted by a blue circle. In general, $s_a$ is different from $s_b$.

Local approximation order can be achieved by using certain higher order of polynomial in the least squares approximation. For example, if quadratic polynomials are used for the local least squares approximation for the LB operator, local truncation error is at least of first order (of second order if super-convergence happens) \cite{liazha13}. However, monotonicity or discrete maximum principle is not preserved by the quadratic approximation. Moreover, the matrix $A$ can be quite ill-conditioned especially when the sampling density of the point cloud is highly non-uniform. Therefore, a more challenging problem is how to construct a discretization so that one can have a computationally efficient algorithm to stably invert the resulting linear system. To address these issues, one approach proposed in  \cite{liazha13} is to enforce the consistency and the diagonal dominance as constraints. Since consistency is guaranteed, which implies that the sum of each row of $A$ is zero, diagonal dominance is enforced if the diagonal element is of different sign with all off-diagonal elements. This means that the local least squares approximation is transformed into a constrained optimization problem. The resulting discretization leads to a M-matrix and hence the discrete maximum principle holds. However, quadratic programming is used for the optimization at each point which increases the computational cost.

\section{A new discretization of the LB operator}
\label{Sec:ModifyApp}

\subsection{The Modified Virtual Grid Difference (MVGD)}

One key idea for the MVGD is to introduce a virtual grid \reminder{aligning} with the local coordinate system of each data point with the grid size adapted to the local sampling density. Instead of using the derivatives of the local least squares approximation at the data point to approximate the corresponding derivatives of the underlying function, we approximate the derivatives of the underlying function by applying an appropriate finite difference method on the local least squares approximation with the following modification: the value of the original function instead of the value of the least squares approximation is used at the center grid which coincides with the data point.

\begin{figure}[!htb]
\begin{center}
\includegraphics[width=6cm]{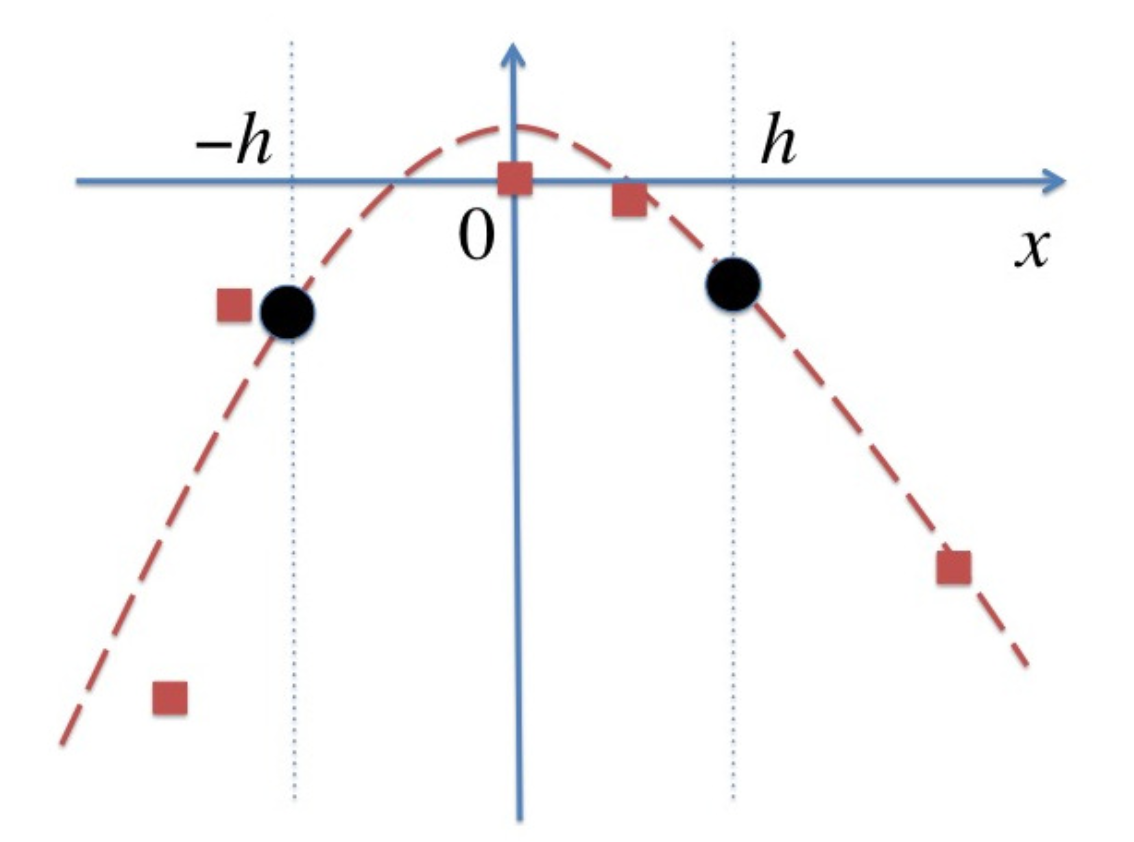}
\end{center}
\caption{Setup for the MVGD. Data from the point cloud ($\p_i$) are plotted in red squares. The least squares approximation to the manifold ($y^{(i)}(x)$) is plotted in the red dashed line. On this parametrized surface, we obtain a least squares approximation to the function $U$, represented by $U^{(i)}(x)$. Then we use the function value at $x=\pm h$ (located using black circles) to approximate the derivatives of $U$ at $x=0$.}
\label{Fig:MVGDSetup}
\end{figure}

Here we first use the following 1D example to illustrate our method. Assume $U$ is a function defined on a 1D curve sampled by a point cloud $S$ with $U_i=U(\p_i)$ and $\p_i\in S$. In a local coordinate with origin at $\p_i$, let $y^{(i)}(x)$ be the least squares approximation of the curve and $U^{(i)}(x)$ be the least squares approximation of $U$ near $\p_i$ in the local coordinate. With a virtual grid in the $x$-dimension centered at $\p_i$ and with a grid size $h$, we define the following MVGD at $\p_i$ to approximate derivatives of $U$,
\begin{eqnarray}
U_{x}(\p_i) &\approx& \tilde{D}_hU(\p_i)\triangleq \frac{U^{(i)}(h)-U^{(i)}(-h)}{2h} \, \mbox{ and } \nonumber\\
U_{xx}(\p_i) &\approx& \tilde{D}^2_hU(\p_i) \triangleq \frac{U^{(i)}(h)-2U_i+U^{(i)}(-h)}{h^2} \, .
\label{eq:centered1}
\end{eqnarray}
Note that we have replaced $U^{(i)}(0)$ by $U_i=U(\p_i)$ for the central finite difference formula for $U^{(i)}_{xx}(0)$. In general, for a least squares approximation, we have $U^{(i)}(0)\ne U_i$ and so $ \tilde{D}^2_hU(\p_i) \ne U^{(i)}_{xx}(0)$. In the case when we use a polynomial with a degree less than 3 in the least squares approximation, we have
\[
U^{(i)}_{x}(0)= \tilde{D}_hU(\p_i)\, .
\]
If the least squares polynomial has a degree less than 4, we have
\[
U^{(i)}_{xx}(0)=\frac{U^{(i)}(h)-2U^{(i)}(0)+U^{(i)}(-h)}{h^2}\ne\tilde{D}^2_hU(\p_i).
\]
For a 1D curve, the LB operator in the local coordinate at $\p_i$ is explicitly given by
\begin{equation}
\Delta_{\Surface}U=\frac{1}{\sqrt{1+y_x^2}}\frac{\partial}{\partial x}\left(\frac{1}{\sqrt{1+y_x^2}}\frac{\partial U}{\partial x}\right)=\frac{U_{xx}}{1+y_x^2}-\frac{y_{x}y_{xx}U_{x}}{\left(1+y_{x}^2\right)^2}\, ,
\label{eq:LB1}
\end{equation}
which leads to the following two possible ways of discretizing the LB operator at $\p_i$.
\begin{enumerate}
\item The central difference
\begin{equation}
\Delta_{\Surface}U\approx \frac{\tilde{D}^2_hU(\p_i)}{1+[y^{(i)}_x(0)]^2}-\frac{y^{(i)}_{x}(0)y^{(i)}_{xx}(0)\tilde{D}_hU(\p_i)}{\left(1+[y^{(i)}_x(0)]^2\right)^2} \, ;
\label{eq:center}
\end{equation}
\item The central difference in the divergence form
\begin{equation}
\frac{ \frac{U^{(i)}(h)}{\sqrt{1+[y^{(i)}_x(\frac{h}{2})]^2}}- \left(\frac{1}{\sqrt{1+[y^{(i)}_x(\frac{h}{2})]^2}} + \frac{1}{\sqrt{1+[y^{(i)}_x(-\frac{h}{2})]^2}}\right) U_i   +\frac{U^{(i)}(-h)}{\sqrt{1+[y^{(i)}_x(-\frac{h}{2})]^2}}}{h^2\sqrt{1+[y^{(i)}_x(0)]^2}} \, .
\label{eq:divergence}
\end{equation}
\end{enumerate}
Both discretizations are easy to implement and have similar performance as will be shown on a simple test in Section \ref{SubSec:LB}. For most of our tests we use the non-divergence form (\ref{eq:center}) in this work.

The MVGD can be easily generalized to high dimensions for the LB operator. To approximate a derivative of $U$ at a point $\p_i$, one can simply apply the standard centered difference to the local least squares approximation $U^{(i)}$ of $U$ on a virtual grid centered at $\p_i$ and replace $U^{(i)}$ at the center grid by $U_i=U(\p_i)$ whenever needed.

The same strategy can be applied to other type of finite difference approximations. For example, to solve the differential equation with an advection term such as the advection diffusion equation, one should use an upwind or a one sided difference. We can made the following modification to the standard one sided difference,
\[
\tilde{D}^+_hU(\p_i)\triangleq \frac{U^{(i)}(h)-U_i}{h} \mbox{ and } \tilde{D}^-_hU(\p_i)\triangleq \frac{U_i-U^{(i)}(-h)}{h}
\]
corresponding to the forward and the background differences, respectively.

The grid scale $h$ of the virtual grid at each data point should be compatible with the local sampling density of the point cloud, which will be justified to some extent below. At a point $\p_i$, the general guideline is that $h$ should not be too large so that $U^{(i)}(h)$ or $U^{(i)}(-h)$ becomes an extrapolation of $U$ in terms of local least squares approximation through its $K$ neighbors $N(i)$. This should be avoided whenever possible since it may cause numerical instability. On the other hand, $h$ should not be too small compared to the spacing of data points near $\p_i$ since it may leads to unnecessary ill-conditioning of the linear system without gaining more accuracy. There are many possible ways to choose $h$ according to this guideline which all produce satisfying and similar numerical results. For example, one could choose $h$ to be the average spacing near $\p_i$ using its KNN. In our numerical experiments, we choose
\begin{equation}
h=\frac{1}{4}(\max\{x_j, \p_j\in N(i)\}-\min\{x_j, \p_j\in N(i)\})
\label{Eqn:DefH}
\end{equation}
for a one dimensional manifold (1D curve), which is one-fourth of the size of the interval that contains the KNN in the local coordinate at $\p_i$. We use
\begin{eqnarray*}
h &=& \min \Big\{\frac{1}{4}\big(\max\{x_j, \p_j\in N(i)\}-\min\{x_j, \p_j\in N(i)\}\big), \\
& & \frac{1}{4}\big(\max\{y_j, \p_j\in N(i)\}-\min\{y_j, \p_j\in N(i)\}\big) \Big\}
\end{eqnarray*}
for a two dimensional manifold.

For point clouds sampled from open surfaces with boundary, various boundary conditions can be incorporated in the MVGD at boundary points. For an open surface $\Sigma$ with the Dirichlet boundary condition $U|_{\partial \Sigma}=f$, we just enforce this condition at all boundary points and use the MVGD for all interior points. For an open surface with the Neumann boundary condition $\partial_{\mathbf{n}} U|_{\partial \Sigma}=f$, where $\mathbf{n}$ is the normal of the boundary in the tangent plane of the surface, we use the following reflection method in the $\mathbf{n}$ direction.

\begin{figure}[!htb]
\centering{
\includegraphics[width=12cm]{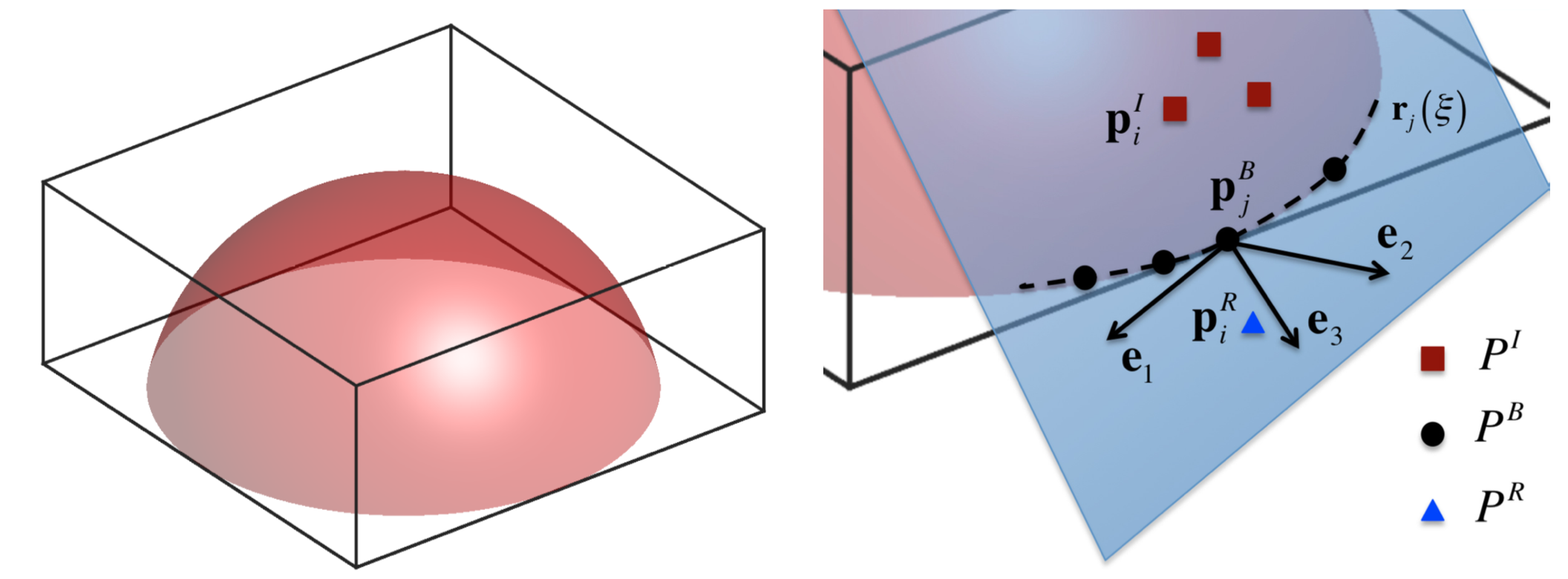}
}
\caption{Notation in the reflection method for imposing the Neumann boundary condition. (Left) The open surface plotted in light red. (Right) The construction of $\p_i^R$ for the interior point $\p_i^I$ which has a boundary point $\p_j^B$ as one of it's KNN.}
\label{Fig:NeumannSetup}
\end{figure}

The key idea of the reflection method is to create an extended layer of ghost points by reflecting those interior points near the boundary with respect to the boundary. The Neumann boundary condition is used to reflect the function values at those interior points near the boundary to their reflected points. Once the reflection is done, the MVGD discretization can be applied at all interior points which results in a linear system that involve all interior points as well as the Neumann boundary condition. In order to implement the reflection for an interior point near the boundary, one needs to construct a local coordinate system at a corresponding boundary point which consists of the normal to the surface, the tangent of the boundary and the normal to the boundary in the tangent plane of the surface. To be more specific mathematically, let $P^I$ and $P^B$ be the sets of interior and boundary points in the point cloud, respectively. For an interior point $\p^I_{i} \in P^I$, we check if any of its KNN is a boundary point. If no, we just use MVGD for the discretization at $\p^I_{i} $ as described before. If yes, we need to create a ghost point as the reflection of $\p^I_{i}$ with respect to the boundary. First, we find the closest boundary point $\p^B_j$, i.e. $j=\mbox{argmin}_k \|\p^I_{i}-\p^B_k\|_2$ and construct a local coordinate system at $\p^B_j$ as follows. Following the procedure proposed in \cite{liazha13}, we use the $k$-nearest boundary neighbors of $\p^B_{j}$ in $P^B$ to construct the boundary curve locally as $\mathbf{r}_j(\xi) =(\xi,\eta_j(\xi),\zeta_j(\xi))$ in some local coordinate system. This provides the tangent of the boundary at $\p^B_{j}$, denoted by $\e_1$. Using the PCA on the $k$-nearest neighbors of $\p^B_{j}$ from the whole data set, one constructs the normal to the surface at $\p^B_{j}$, denoted by $\e_2$. Then the normal to the boundary at $\p^B_{j}$ in the tangent plane of the surface is then given by $\e_3=\e_1\times \e_2$. Assuming $\p^I_{i}$ has coordinates $(x,y,z)$ in this local coordinate system, its reflection point $\p^R_{i}$ has coordinates $(x,y,-z)$. The setup is summarized and  is plotted in Figure \ref{Fig:NeumannSetup}.

Once the reflection point is determined, we impose the Neumann boundary condition $\partial_{\mathbf{n}} U|_{\partial \Sigma}=f$ by assigning the function value at the reflected point by
$$
U(\p^R_{i})=U(\p^I_{i})+\left\|\p^I_{i}-\p^R_{i}\right\| \, \cdot f(\p^B_{j}) \, .
$$
Once we obtain all reflections of all interior points near the boundary with respect to the boundary, we assign their function values by imposing the boundary condition and we use the MVGD to discretize the differential operator at each interior point $\p^I_i\in P^I$. Note that the KNN of $\p^I_i$ consist of points in $P^I \cup P^R$, where $P^R$ denotes the set of all reflected points. Hence a linear system involving all interior points is formed with the Neumann boundary condition implicitly incorporated. Combining the treatments for the Dirichlet and the Neumann boundary conditions, one can also easily deal with the general Robin type of boundary condition.

Finally, for the full linear system resulting from the above discretization for the LB operator at each point, we find that algebraic multi-grid (AMG) is a very effective solver.

\subsection{Some motivations and analysis of the MVGD}
\label{SubSec:OneDEx}
For simplicity, we use the Laplacian in $\mathbb{R}^1$ to shed some insights into the new discretization in terms of both accuracy and stability.
\subsubsection{Accuracy}
In terms of local truncation error, from the definitions in (\ref{eq:centered1}), we have
\begin{equation}
 \tilde{D}^2_hU(\p_i) =  \frac{U^{(i)}(h)-2U^{(i)}(0)+U^{(i)}(-h)}{h^2}+\frac{2(U^{(i)}(0)-U_i)}{h^2} \, .
 \label{eq:accuracy}
\end{equation}
And so, one can see that when $h$ is compatible with the local spacing of data points and quadratic polynomial is used for least squares approximation, the two terms in the above expression for local truncation error is of the same order. They are at least of $O(h)$ and can be of $O(h^2)$ if the distribution of data points has some symmetry so that super-convergence occurs (see \cite{liazha13}).

\subsubsection{Stability}
As discussed before, an important and more difficult issue is for a discretization is to result a better behaved linear system that can be solved stably and efficiently. Since we are discretizing the LB operator directly on a point cloud without a global mesh or parametrization, it is rather hard, if not impossible, to design a discretization so that the resulting linear system has a symmetric positive definite matrix. However, one can still hope that the linear system has a M-matrix due to the use of least squares approximation instead of exact interpolation. Discretization based on the least squares approximation provides (i) more robustness with respect to noise or almost degeneracy of point distribution, e.g., very close neighbors when sampling is highly non-uniform, and (ii) flexibility to satisfy both required accuracy and extra constraints. For example, a constrained optimization approach was proposed in \cite{liazha13} to design such a discretization. However, quadratic programming has to be used for the optimization at each point. Here we show that the proposed MVGD for the LB operator can also achieve this goal.

Since our discretization of the LB operator is consistent,  the sum of each row of the matrix of the discretized linear system is zero. If the matrix is diagonal dominant, then it is a M-matrix. Below we construct an explicit example for the Laplacian in $\mathbb{R}^1$ to show that by using the value of the original function instead of the value of the least squares approximation at the center grid in the standard central finite difference scheme, diagonal dominance can be achieved for the linear system when the grid size is comparable to the local spacing of the point cloud.

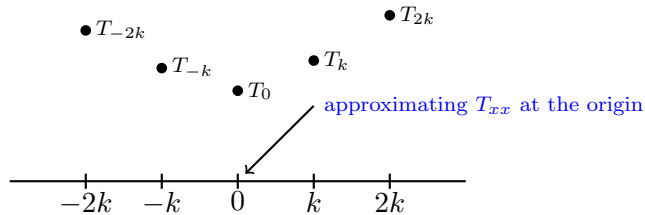
\begin{figure}[!htb]
\centering{
\begin{tikzpicture}
\draw [thick,-] (-3,0) -- (3,0);
\draw [thick,-] (-2,-0.1) -- (-2,0.1);
\draw [thick,-] (-1,-0.1) -- (-1,0.1);
\draw [thick,-] (0,-0.1) -- (0,0.1);
\draw [thick,-] (1,-0.1) -- (1,0.1);
\draw [thick,-] (2,-0.1) -- (2,0.1);
\fill (-2,2) circle[radius=2pt];
\fill (-1,1.5)  circle[radius=2pt];
\fill (0,1.2)  circle[radius=2pt];
\fill (1,1.6)  circle[radius=2pt];
\fill (2,2.2)  circle[radius=2pt];
\node [below] at (-2,0) {$-2k$};
\node [below] at (-1,0) {$-k$};
\node [below] at (0,0) {0};
\node [below] at (1,0) {$k$};
\node [below] at (2,0) {$2k$};
\node [right] at (-2,2) {\scriptsize$T_{-2k}$};
\node [right] at (-1,1.5) {\scriptsize$T_{-k}$};
\node [right] at (0,1.2) {\scriptsize$T_{0}$};
\node [right] at (1,1.6) {\scriptsize$T_{k}$};
\node [right] at (2,2.2) {\scriptsize$T_{2k}$};
\draw [thick,<-] (0.1,0.1) -- (1,1);
\node [right,blue] at (1,1) {\scriptsize approximating $T_{xx}$ at the origin};
\end{tikzpicture}
\caption{A 1D example.}
\label{Fig:Analysis}
}
\end{figure}

Consider that the discretization of the Laplacian of a function $T$ on the straight line at the point $x=0$ as shown in Fig \ref{Fig:Analysis}, and the function $T$ is sampled at regular grid points with grid size $k$.

Assume $T^{(0)}(x)=a_0+a_1x+a_2x^2$ is the least squares quadratic approximation through the function values at 5 points, i.e. $T_j=T(jh)$ for $j=0, \pm1, \pm2$.  From the least squares approximation, we have
$$
 \left( \begin{array}{ccc}
 5 & 0 & 10k^2\\
 0 & 10k^2 & 0\\
 10k^2 & 0 & 34k^4
 \end{array} \right)
\left( \begin{array}{c}
a_0\\
a_1\\
a_2
\end{array} \right) =
\left( \begin{array}{ccccc}
1 & 1 & 1 & 1 & 1\\
-2k & -k & 0 & k & 2k\\
4k^2 & k^2 & 0 & k^2 & 4k^2
\end{array} \right)
\left( \begin{array}{c}
T_{-2k}\\
T_{-k}\\
T_{0}\\
T_{k}\\
T_{2k}
\end{array} \right) \, ,
$$
and so
$$
\left( \begin{array}{c}
a_0\\
a_1\\
a_2
\end{array} \right) = \frac{1}{70k^2}
\left( \begin{array}{ccccc}
-6k^2 & 24k^2 & 34k^2 & 24k^2 & -6k^2\\
-14k & -7k & 0 & 7k & 14k\\
10 & -5 & -10 & -5 & 10
\end{array} \right)
\left( \begin{array}{c}
T_{-2k}\\
T_{-k}\\
T_{0}\\
T_{k}\\
T_{2k}
\end{array} \right) \, .
$$
If we simply use the derivative of the least squares approximation at the origin, it leads to the following discretization of $T_{xx}$
$$
T_{xx}(0) \approx T^{(0)}_{xx}(0) =2a_2=\frac{1}{7k^2}(2,-1,-2,-1,2) {\mathbf {T}}
$$
where ${\mathbf {T}}=(T_{-2k},T_{-k},T_{0},T_{k},T_{2k})^T$. Although it is consistent with the second order local truncation error, the resulting matrix is not diagonal dominant since off-diagonal elements have mixed sign. This means that the linear system after discretization does not give a M-matrix and this can cause numerical instability and might lead to non-convergence for many iterative methods.

Now for the proposed MVGD discretization, we have
\begin{eqnarray*}
T_{xx}(0) &\approx& \tilde{D}^2_hT^{(0)}(0) =  \frac{T^{(0)}(h)-2T_0+T^{(0)}(-h)}{h^2}
=\frac{2a_0-2T_0}{h^2}+2a_2 \\
&=& \frac{2}{70 \, h^2} \left( -6 + \frac{10h^2}{k^2}, 24 - \frac{5h^2}{k^2},
-36 - \frac{10h^2}{k^2}, 24 - \frac{5h^2}{k^2}, -6 + \frac{10h^2}{k^2} \right) \mathbf{T} \, .
\end{eqnarray*}
As long as
\begin{equation}
\frac{3}{5} \le \left(\frac{h}{k} \right)^2 \le \frac{24}{5},
\label{eq:condition}
\end{equation}
all off-diagonal elements have different sign from diagonal element which makes the discretized linear system having a M-matrix. Since $\max x_i = 2k$ and $\min x_i = -2k$ in this simple example, our choice of $h$ according to (\ref{Eqn:DefH}) implies $h=k$ and so the ratio $(h/k)^2$ satisfies above conditions.

From this simple yet explicit example, we see that by replacing $T_0$ by $T^{(0)}(0)$ at the center grid in the modified central finite difference, we have introduced to the coefficient $a_2$ in the new discretization an extra term, $(a_0-T_0)/h^2$. It can be easily checked that the magnitude of the newly added term is of order $O(h^2)$ which is the same as $|a_2-T_{xx}(0)|$ when $h\sim k$, which is a special case of equation \eqref{eq:accuracy} when a quadratic polynomial is used for the least squares approximation and data points are regularly distributed around 0. However, the key observation is that this modification enhances the diagonal dominance. In this special case, for a given $k$, diagonal dominance is guaranteed for any $h$ satisfying \eqref{eq:condition}.
For more general cases, an explicit formula is more complicated since it depends on both the distribution of neighboring points near 0 and the number of points used in the least squares approximation. Nevertheless, we will demonstrate numerically that the new discretization of LB operator on point clouds will indeed lead to a much better conditioned matrix in the following section.

Finally, in practice, at each point one could also determine the virtual grid spacing $h$ by the following dynamical approach. We can first determine explicitly all the coefficients in the discretization for the LB operator in terms of $h$, and then we choose the value of $h$ at each sampling point to optimize the diagonal dominant property.


\section{Numerical examples}
\label{Sec:Examples}

In this section, we give numerical results from our proposed MVGD discretization of the LB operator on point clouds. We are going to demonstrate that the new discretization indeed provides better conditioning, stability and numerical accuracy when compared to the simple least squares approach (without using the constraint optimization) as described in Section \ref{Sec:Background} as proposed in \cite{liazha13}. For all numerical examples, the linear system corresponding to the MVGD discretization is solved by the AGM. While the linear system corresponding to the simple least squares approach cannot be typically solved by most common efficient iterative methods and is solved by GMRES in our tests.

\subsection{Conditioning of the discretized system}
\label{sec:conditioning}
Let $A$ be the matrix corresponding to the discretization of the LB operator on a given point cloud. We decompose it as $A=M+N$ according to Gauss-Seidel iteration, i.e., $N$ is upper triangular portion of $A$. We compute the eigenvalues of matrix $M^{-1}N$ for a few point clouds in 2D and 3D.

\begin{figure}[!htb]
\begin{center}
\includegraphics[width=12cm]{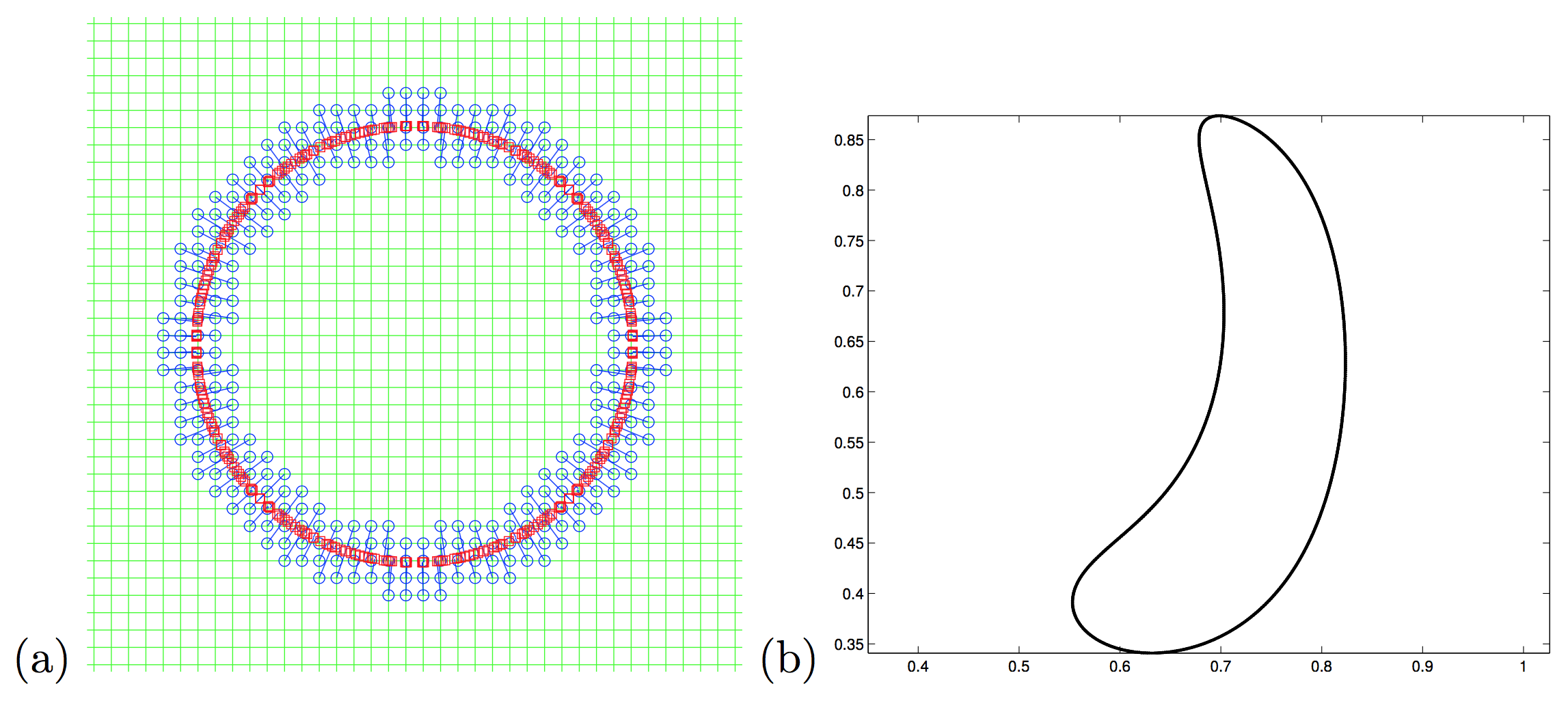}
\end{center}
\caption{(a) A circle sampled by the GBPM. (b) Another closed curve used in the test.} \label{Fig:twocurves}
\end{figure}

\begin{figure}[!htb]
\begin{center}
\includegraphics[width=12cm]{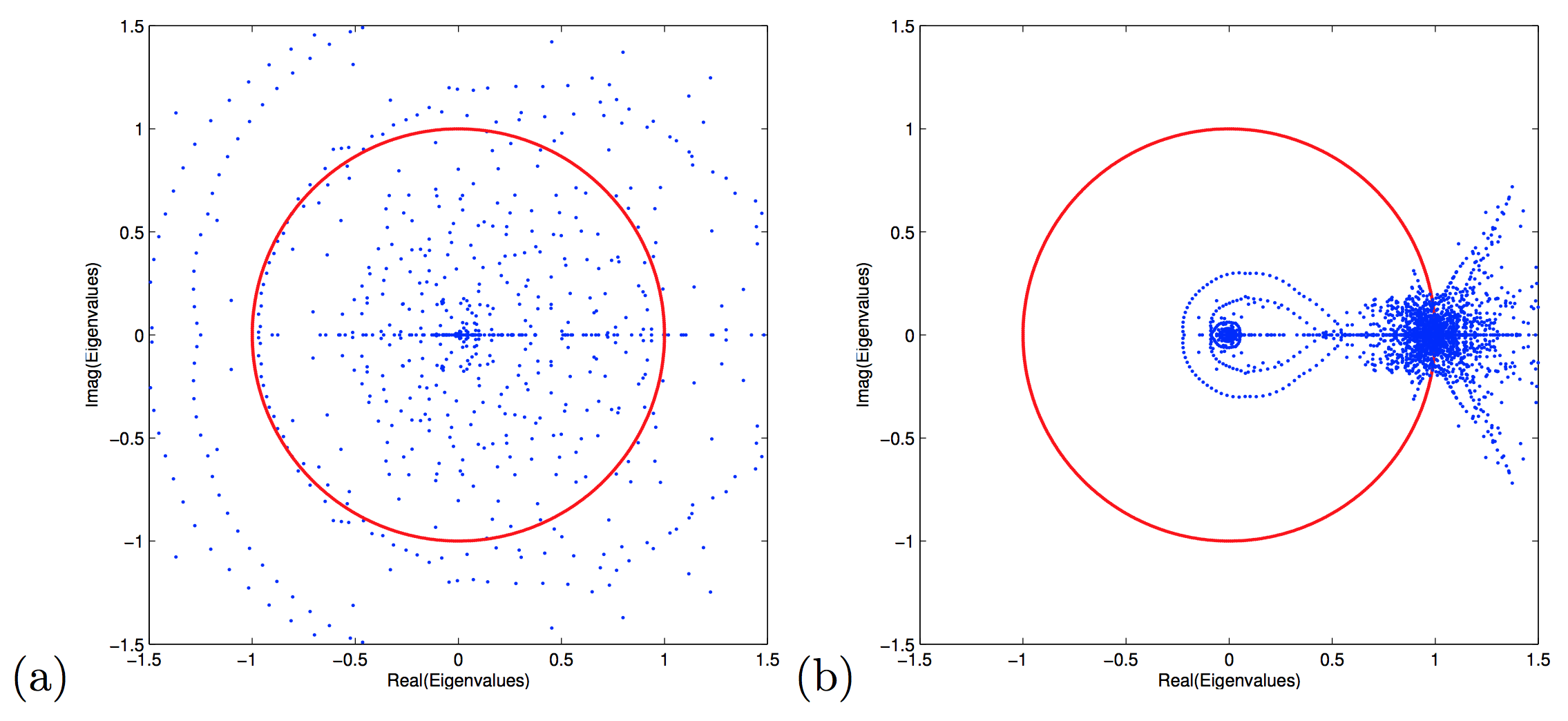}
\end{center}
\caption{Eigenvalues of the Gauss-Seidel iterative matrix for LS on (a) a circle and (b) a closed curve as shown in Fig \ref{Fig:twocurves} (b).
}
\label{Fig:Pre_LapSpectrum}
\end{figure}

\begin{figure}[!htb]
\begin{center}
\includegraphics[width=12cm]{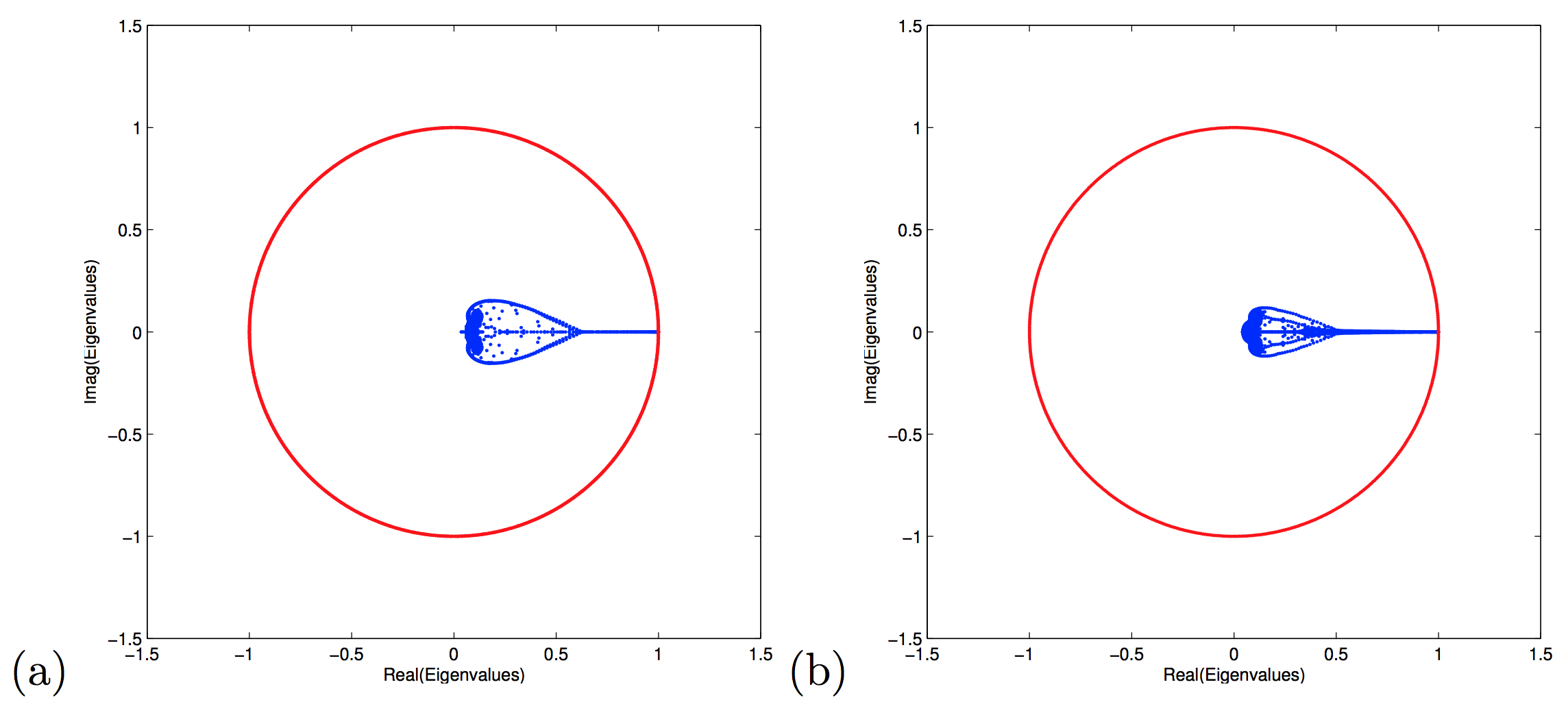}
\end{center}
\caption{Eigenvalues of the iterative matrix for MVGD on (a) a circle and (b) a closed curve as shown in Fig \ref{Fig:twocurves} (b).}
\label{Fig:LapSpectrum}
\end{figure}

In our first example, we have sampled the two curves in Fig \ref{Fig:twocurves} by the GBPM, where an interface is sampled by closest points to those underlying mesh points in the vicinity of the interface. A typical point cloud sampling an interface by the GBPM is illustrated in Fig \ref{Fig:twocurves} (a), where we plot the underlying mesh in solid line, all active grids near the interface (a circle) using little blue circles and their associated closest points on the interface using little red squares. The correspondence between each pair is shown by a solid line link. The point cloud resulting from the GBPM can be quite non-uniform because two closest points to two mesh points can be very close or the same. The eigenvalues of matrix $M^{-1}N$ corresponding to the discretization using the LS approach for these two point clouds are plotted in Fig \ref{Fig:Pre_LapSpectrum}. The unit circle in the complex plane is plotted in red. As we can see clearly, the magnitude of most of these eigenvalues (blue dots) are larger than one and there is no guarantee of convergence if Gauss-Seidel iteration is used to solve the linear system. Fig \ref{Fig:LapSpectrum} shows the  eigenvalues corresponding to the MVGD discretization. The magnitude of all eigenvalues are now less than 1. Hence, even the simple Gauss-Seidel iteration can be used to solve the linear system.

\subsection{Solving the LB equation on point clouds}
\label{SubSec:LB}

\begin{figure}[!htb]
\begin{center}
\includegraphics[width=12cm]{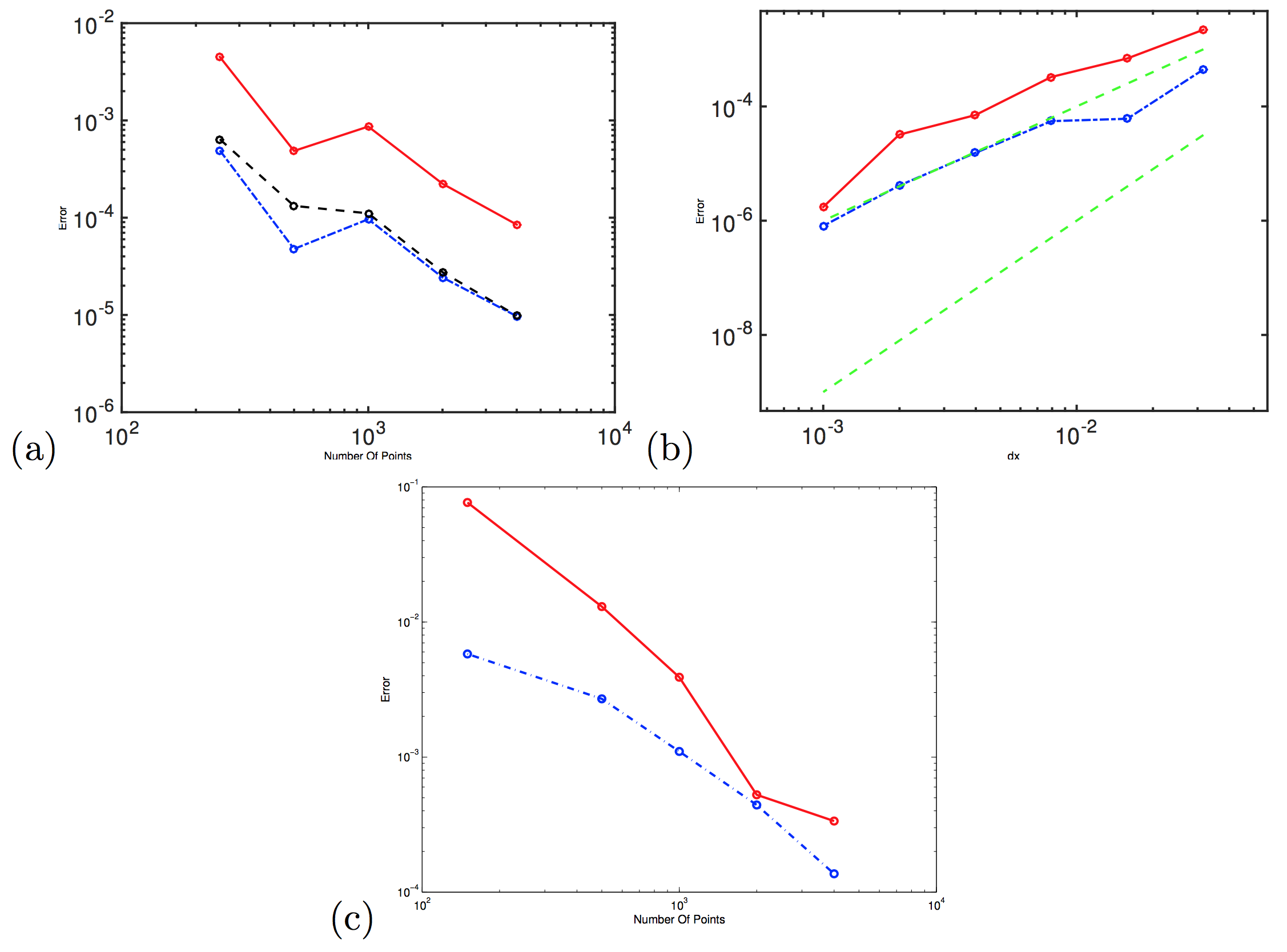}
\caption{Error in the numerical solution to the LB equation on point clouds sampled from a circle using (a) the uniform  sampling, (b) the GBPM sampling \reminder{and (c) the non-uniform sampling}. Red line: MLS. Blue dash-dotted line: MVGD in the non-divergence form. Black dash line: MVGD in the divergence form. Green lines: reference curves of $O(dx^2)$ and $O(dx^3)$, respectively.}
\label{Fig:AMG_Uniform}
\end{center}
\end{figure}

We first consider a simple example, point clouds sampled from a unit circle centered at the origin. On the unit circle, let $U(\theta)=-2\sin \theta \cos \theta$ which solves $-\Delta_{\Surface} U=-8 \sin \theta \cos \theta$. We solve this LB equation on three different point cloud data sets sampled on the circle. One is evenly distributed, the second one is sampled by the GBPM \reminder{and the third one is a non-uniformly sample obtained by randomly picking points on the circle.} 

\reminder{Fig \ref{Fig:AMG_Uniform} shows the $L^{\infty}$ errors in the numerical solution using these three different type of sampling methods. In Fig \ref{Fig:AMG_Uniform} (a) and (c), we plot the error in the solutions versus the number of sampling points. In Fig \ref{Fig:AMG_Uniform} (b), we consider the error versus the grid spacing $dx$ in the underlying uniform mesh. As a comparison, we have also implemented the MLS approach as discussed in \cite{liazha13}. The weights in the method are chosen to be $w(d)=1$ if $d=0$ and $w(d)=1/K$ if $d \neq 0$, where $K$ is the number of neighboring points, which is reported to be one of the best weighting functions.} Even though all methods seem to give second order convergence, the proposed discretization in this paper can achieve a better accuracy in solving the LB equation. More importantly, since the linear system resulted from the MVGD discretization can be solved by any efficient iterative solver, such as AMG, it is computationally much more efficient. We have implemented the MVGD for the LB operator using both central difference (\ref{eq:center}) and the divergence form (\ref{eq:divergence}) and show results in Fig \ref{Fig:AMG_Uniform} (a). The two are comparable and both are significant better than that by using direct MLS discretization.

\begin{figure}[!htb]
\includegraphics[width=12cm]{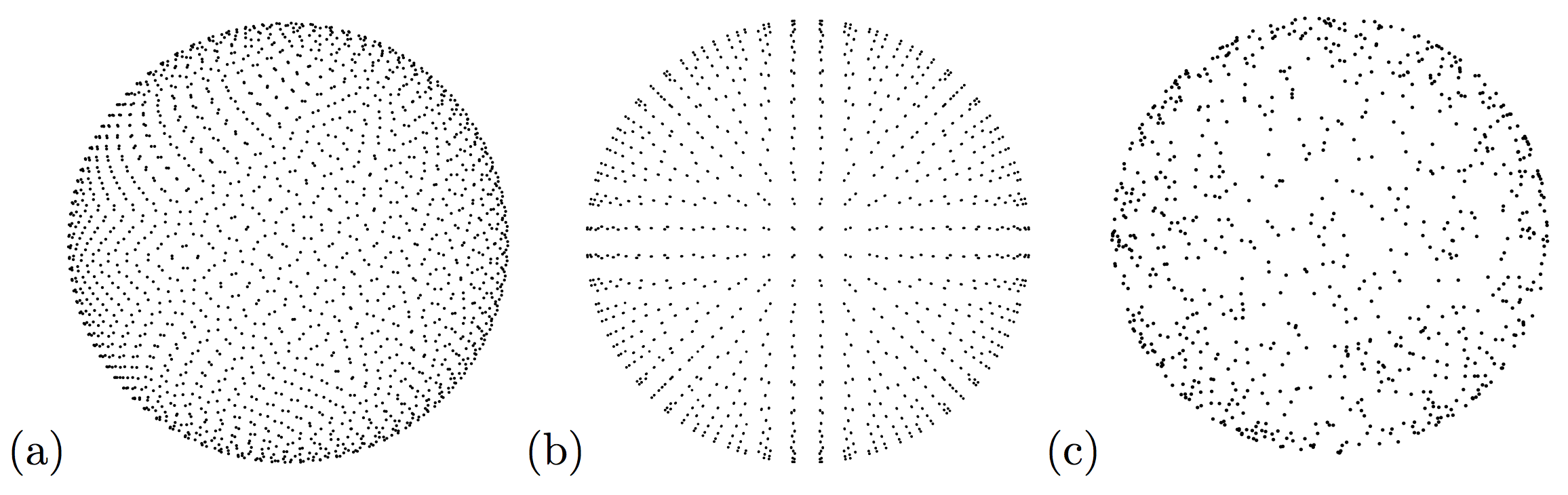}
\caption{Different samplings of a unit sphere. (a) The Fibonacci sampling, (b) the GBPM sampling, \reminder{ and (c) a non-uniform sampling.}}
\label{Fig:sampling}
\end{figure}

We carry out a similar test on a sphere. Let $\p$ be a point on any given surface $\Surface$, we have
\begin{equation}
-\Delta_{\Surface} \p=2H(\p) \label{eqn-lapsphere}
\end{equation}
where $H(\p)$ is the mean curvature vector on $\Surface$ \cite{Willmore93}. If $\Surface$ is the unit sphere centered at the origin, the mean curvature is one and the above equation is reduced to $-\Delta_{\Surface} \p=2\p$. To test the convergence of the algorithm, we look at only the $x$-coordinate in $\p$, i.e.  we let the right hand side of (\ref{eqn-lapsphere}) be $2x$ and therefore the exact solution to the equation is simply $x$. Once again, we consider three different point clouds. One is the \textit{uniform} sampling given by the Fibonacci sampling \cite{gon10}, \reminder{the second one is sampled by the GBPM and the third one is non-uniform sampled by randomly picking points on a sphere}. These three different samplings of the sphere are shown in Fig \ref{Fig:sampling}.

\begin{figure}[!htb]
\begin{center}
\includegraphics[width=12cm]{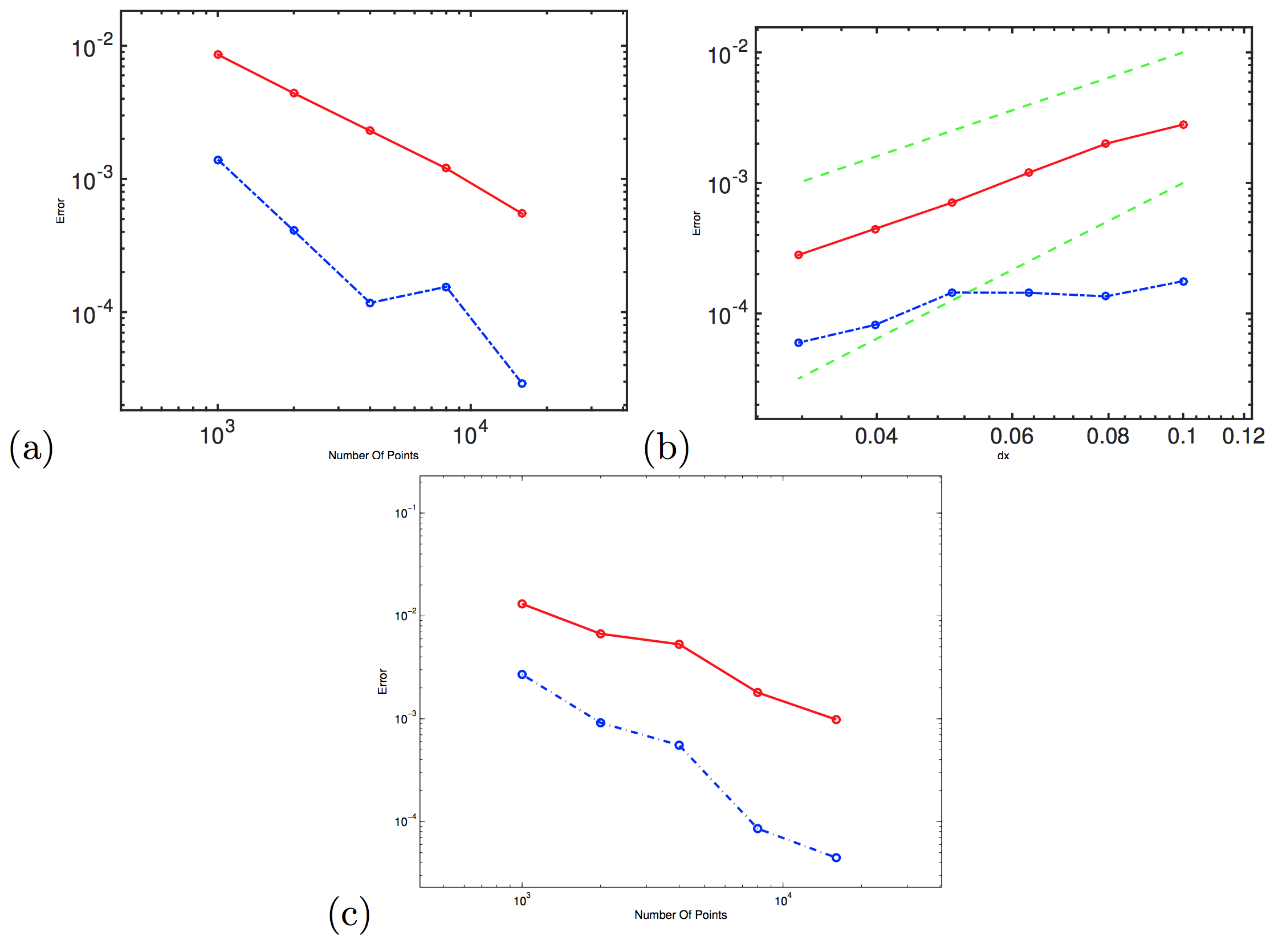}
\caption{\reminder{Error in the numerical solution to the LB equation on a sphere sampled by (a) the Fibonacci sampling, (b) the GBPM sampling, and (c) the non-uniform Sampling}. Red line: MLS. Blue dash-dotted line: MVGD. Green lines: reference curves of $O(dx^2)$ and $O(dx^3)$, respectively.}
\label{Fig:AMG_GBPM}
\end{center}
\end{figure}

Fig \ref{Fig:AMG_GBPM} (a) shows the $L^{\infty}$ error versus the number of points on a uniformly sampled sphere, and (b) shows the corresponding $L^{\infty}$ error versus the grid spacing $dx$ on a sphere sampled by the GBPM \reminder{and (c) shows the $L^{\infty}$ error versus the number of points on a non-uniformly sampled sphere}. We observe similar convergence behaviors as in the circle case. Both methods converge in approximately second order. AMG works beautifully for the linear system discretized by the MVGD. Table \ref{table:TableCPUTime} presents the CPU time (in second) for solving the LB equation on a sphere sampled by the uniform sampling. It shows that the CPU time is approximately linear in the number of sampling points.

\begin{table}[!htb]
\begin{center}
\begin{tabular}{|c||c|c|c|c|c|c|}
\hline
Sample size & 500 & 1000 & 2000 & 4000 & 8000 & 16000 \\
\hline
CPU time (s) & 0.0239 & 0.0512 & 0.133 & 0.261 & 0.365 & 0.845\\
\hline
\end{tabular}
\caption{The CPU time for solving LB equation discretized by the MVGD using the AMG. Sample size represents the number of points on a uniformly sampled sphere. The CPU time is approximately linearly proportional to the sample size.}
\label{table:TableCPUTime}
\end{center}
\end{table}

\subsection{Eigenvalues and eigenfunctions of the  LB operator}
\label{LBEig}

In this section, we solve the eigenvalue problem for the LB operator on manifolds, i.e we determine constants $\lambda$ and corresponding functions $\nu$ such that
$$
-\Delta_{\Sigma} \nu = \lambda \nu
$$
on a given closed surface $\Sigma$. Numerically, we discretize the LB operator and determine the eigenvalues $\lambda$ and their corresponding eigenvectors $\nu$ such that $A\nu=\lambda \nu$ using the function \textsf{eigs} in \textsf{MATLAB}. For the LB eigenvalue problem on a unit sphere, the exact eigenvalues and their corresponding eigenvectors are explicitly known. The $n$-th eigenvalue is given by $\lambda_{n} = n(n+1)$ with multiplicity $(2n+1)$ and the associated eigenfunction is given by the spherical harmonics.

\begin{table}[!htb]
\begin{center}
\scalebox{0.775}{
\begin{tabular}{|c||c|c|c|c||c|c|c|c||c|c|c|c|}
\hline
&\multicolumn{4}{|c||}{Uniform} & \multicolumn{4}{|c||}{GBPM} & \multicolumn{4}{|c|}{Random}\\
\hline
samples & 2000 & 4000 & 8000 & 16000 & 2007 & 4011 & 8048 & 16038 & 2000 & 4000 & 8000 & 16000\\
\hline
\multicolumn{13}{|c|}{MVGD}\\
\hline
$\lambda_4=20$ & 0.48 & 0.15 & 0.12 & 0.050 & 1.2 & 0.76 & 0.35 & 0.21 & 0.98 & 0.41 & 0.19 & 0.069\\
\hline
$\lambda_8=72$ & 1.2 & 0.77 & 0.34 & 0.21 & 1.9 & 0.90 & 0.76 & 0.51 & 3.12 & 1.38 & 0.63 & 0.26\\
\hline
\multicolumn{13}{|c|}{MLS \cite{liazha13}}\\
\hline
$\lambda_4=20$ & 2.1 & 1.1 & 0.54 & 0.27 & 1.5 & 0.74 & 0.37 & 0.18 & 3.91 & 2.1 & 0.70 & 0.41\\
\hline
$\lambda_8=72$ & 26 & 5.0 & 2.5 & 1.2 & 6.8 & 3.5 & 1.7 & 0.85 & 27.8 & 10.23 & 3.41 & 1.97\\
\hline
\multicolumn{13}{|c|}{FEM \cite{mdsb03}}\\
\hline
$\lambda_4=20$ & 0.69 & 0.31 & 0.13 & 0.063 & \multicolumn{8}{|c|}{}\\
\hline
$\lambda_8=72$ & 2.7 & 1.2 & 0.62 & 0.32 & \multicolumn{8}{|c|}{}\\
\hline
\multicolumn{13}{|c|}{LMM \cite{lailiazha13}}\\
\hline
$\lambda_4=20$ & 0.89 & 0.45 & 0.21 & 0.10 & 2.01 & 1.13 & 0.67 & 0.53 & 3.21 & 2.08 & 1.17 & 0.76\\
\hline
$\lambda_8=72$ & 2.78 & 1.49 & 0.71 & 0.32 & 6.71 & 3.12 & 1.02 & 0.65 & 9.76 & 4.78 & 2.20 & 1.98\\
\hline
\end{tabular}
}
\caption{\reminder{$E^2_n$ ($\times10^{-2}$) on different samplings of a sphere. We compare our proposed MVGD approach with the moving least squares method (MLS) \cite{liazha13}, the finite element method (FEM) \cite{mdsb03} and also the local mesh method (LMM) \cite{lailiazha13}.}}
\label{table:EigTableProposed2}
\end{center}
\end{table}

\begin{table}[!htb]
\begin{center}
\scalebox{0.775}{
\begin{tabular}{|c||c|c|c|c||c|c|c|c||c|c|c|c|}
\hline
&\multicolumn{4}{|c||}{Uniform} & \multicolumn{4}{|c||}{GBPM} & \multicolumn{4}{|c|}{Random}\\
\hline
samples & 2000 & 4000 & 8000 & 16000 & 2007 & 4011 & 8048 & 16038 & 2000 & 4000 & 8000 & 16000\\
\hline
\multicolumn{13}{|c|}{MVGD}\\
\hline
$\lambda_4=20$ & 0.72 & 0.27 & 0.25 & 0.0765 & 1.57 & 1.05 & 0.59 & 0.38 & 1.69 & 0.74 & 0.30 & 0.13\\
\hline
$\lambda_8=72$ & 1.76 & 1.46 & 0.74 & 0.38 & 3.33 & 1.68 & 1.67 & 1.22 & 4.89 & 2.78 & 1.33 & 0.47\\
\hline
\multicolumn{13}{|c|}{MLS \cite{liazha13}}\\
\hline
$\lambda_4=20$ & 2.49 & 1.13 & 0.57 & 0.29 & 1.56 & 0.8 & 0.42 & 0.23 & 5.11 & 2.57 & 0.89 & 0.65\\
\hline
$\lambda_8=72$ & 28.25 & 5.53 & 2.63 & 1.36 & 7.3 & 3.73 & 1.9 & 0.98 & 42.53 & 18.79 & 6.71 & 3.01\\
\hline
\multicolumn{13}{|c|}{FEM \cite{mdsb03}}\\
\hline
$\lambda_4=20$ & 0.82 & 0.41 & 0.21 & 0.11 & \multicolumn{8}{|c|}{}\\
\hline
$\lambda_8=72$ & 3.51 & 1.70 & 0.85 & 0.42 & \multicolumn{8}{|c|}{}\\
\hline
\multicolumn{13}{|c|}{LMM \cite{lailiazha13}}\\
\hline
$\lambda_4=20$ & 1.03 & 0.50 & 0.25 & 0.13 & 2.57 & 1.37 & 0.86 & 0.65 & 3.67 & 2.42 & 1.37 & 0.91\\
\hline
$\lambda_8=72$ & 3.57 & 1.67 & 0.87 & 0.43 & 8.63 & 4.12 & 1.36 & 0.97 & 10.76 & 5.97 & 2.97 & 2.16\\
\hline
\end{tabular}
}
\caption{\reminder{$E^{\infty}_n$ ($\times10^{-2}$) on different samplings of a sphere. We compare our proposed MVGD approach with the moving least squares method (MLS) \cite{liazha13}, the finite element method (FEM) \cite{mdsb03} and also the local mesh method (LMM) \cite{lailiazha13}.}}
\label{table:EigTableProposedInf}
\end{center}
\end{table}

\reminder{Again we have tested on three different point clouds of a sphere including the uniform sample, the sampling by the GBPM and also the non-uniform random sampling}. For a particular eigenvalue $\lambda_{n}$, we define the following $L^2$- and $L^{\infty}$-norm error by
$$
E^{2}_n=\sqrt{\frac{1}{2n+1}\sum_i \left(\frac{\lambda_{n,i}-\lambda_n}{\lambda_n} \right)^2} \, \mbox{ and } \, E^{\infty}_n=\max_i \left| \frac{\lambda_{n,i}-\lambda_n}{\lambda_n} \right| \, ,
$$
where $\lambda_{n,i}$ are the eigenvalues computed from the the discretized LB matrix $A$ for eigenvalue, and $i=1,2,...,2n+1$. Table \ref{table:EigTableProposed2} and Table \ref{table:EigTableProposedInf} \reminder{compare our proposed discretization with the MLS method \cite{liazha13}, the finite element method (FEM) as proposed in \cite{mdsb03} and also the local mesh method (LMM) developed in \cite{lailiazha13}}. In each of these tables, we look at the errors in the fourth and the eighth eigenvalues (i.e. $\lambda_4=20$ and $\lambda_8=72$).

For the Fibonacci uniform sampling \reminder{and the non-uniform random sampling}, the MVGD based approach gives more accurate solutions in all sampling densities. When the sampling density of the point cloud \reminder{is obtained by the GBPM, the accuracy is comparable. For the finite element method in \cite{mdsb03}, we have omitted those tests when the sampling is non-uniform as in the GBPM or in the random sampling since a global mesh cannot be obtained.} Again, the more important point is that the linear system corresponding to the new discretization is better conditioned and can be solved much more efficiently by using off-the-shelf fast solver such as the AMG.

Note, however, that since the matrix is not symmetric, there is no guarantee that all obtained eigenvalues are real. In fact, we have indeed obtained some complex eigenvalue pairs in the solutions. For example, considering the sphere represented by the uniform sampling, we found one complex pair among the first 81 eigenvalues for the data with 2000 samples, and two complex pairs for each of the samplings with 4000, 8000 and 16000 sampling points. Nevertheless, the largest relative magnitude of the complex parts are given by $1.24\times10^{-4}$, $1.29\times10^{-5}$, $1.05\times10^{-6}$, and $1.71\times10^{-5}$ for four different sampling densities, respectively. When the surface is sampled by the GBPM, we have a slightly larger number of complex pairs. Among the first 81 eigenvalues for the same number of sampling points as in Table \ref{table:EigTableProposed2}, we found 9, 4, 7 and 3 complex eigen-pairs with the largest relative magnitudes given by $2.14\times10^{-3}$, $3.22\times10^{-3}$, $1.34\times10^{-3}$, and $7.56\times10^{-4}$ for four different sampling densities, respectively. When we determine the error in the eigenvalues, we simply truncate the imaginary part of these complex eigen-pairs.

\begin{figure}[!htb]
\centering{
\includegraphics[width=12cm]{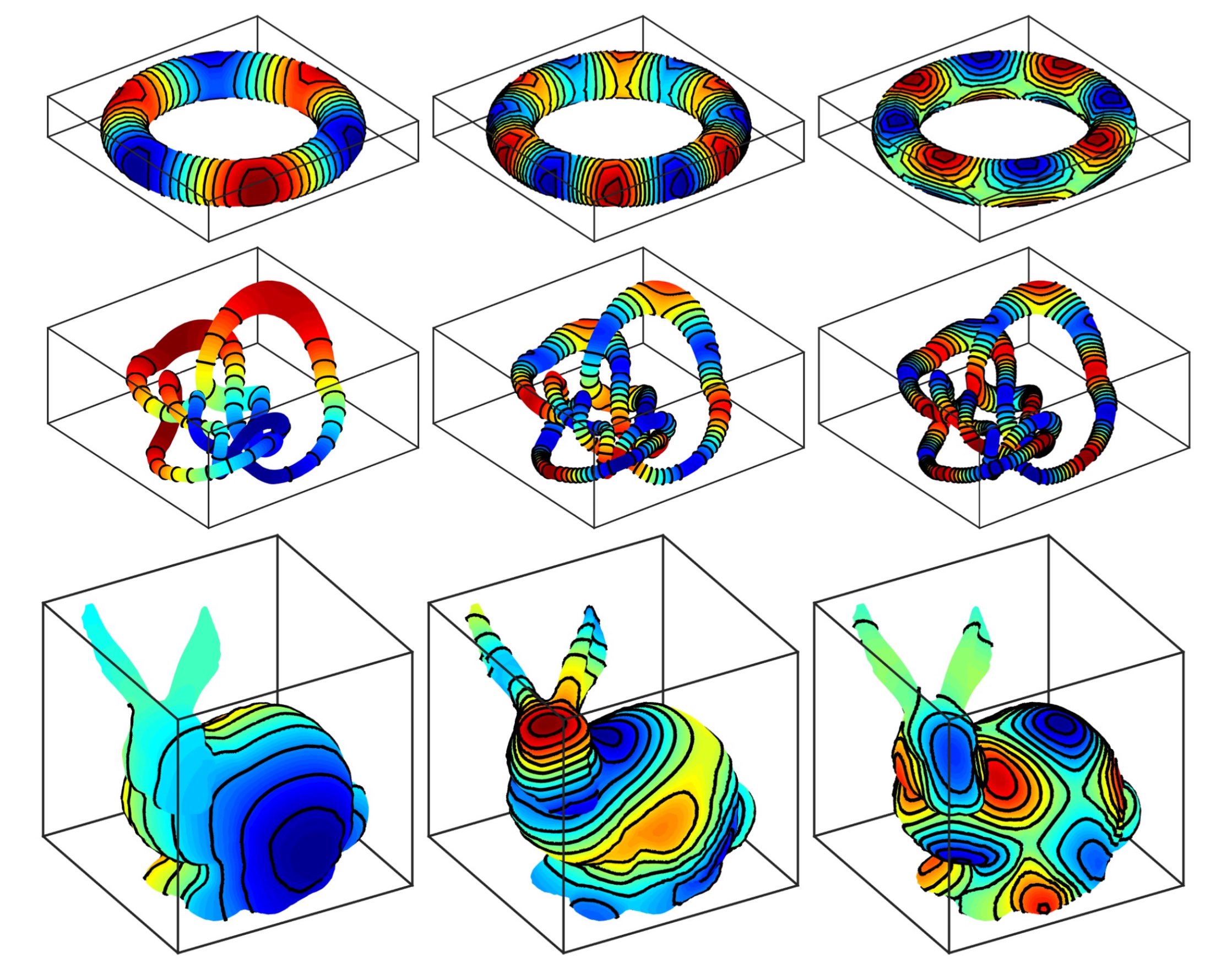}
}
\caption{The sixth, the twenty-second and the thirtieth eigenfunctions of the LB operator on different point clouds.}
\label{Fig:Eig_PointCloud}
\end{figure}

We also tested our new method to compute the eigenvalue problem of the LB operator on several publicly available point clouds including a torus (768 points), a knot (36898 points) and also the Stanford bunny (28006 points). We have computed the first 30 smallest eigenvalues (in magnitude) and their corresponding eigenfunctions. The sixth, the twenty-second and the thirtieth eigenfunctions from the \textsf{MATLAB} function \textsf{eigs} are shown in Fig \ref{Fig:Eig_PointCloud}.

\begin{figure}[!htb]
\centering{
\includegraphics[width=12cm]{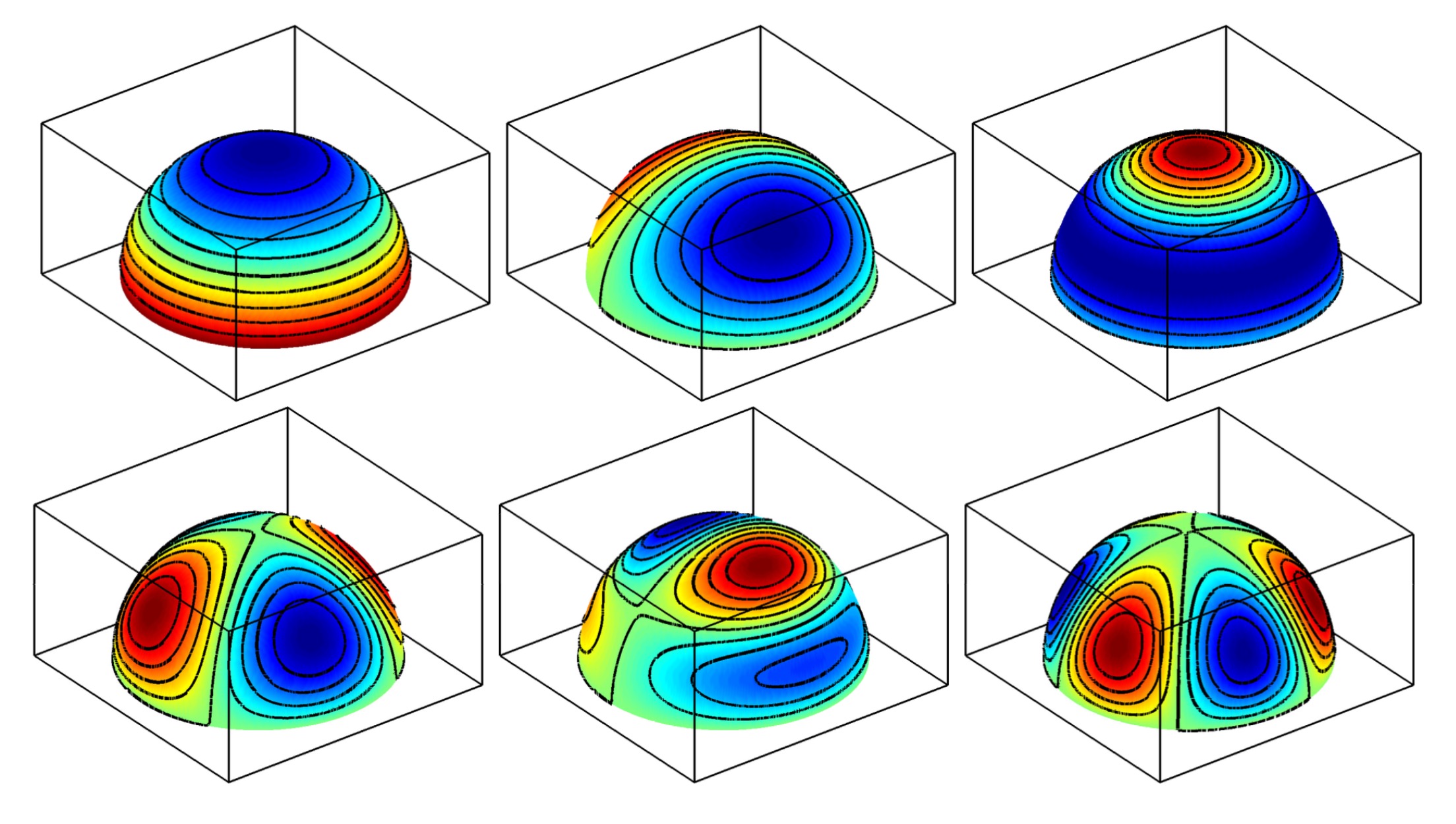}
}
\caption{First few eigenfunctions of the LB operator with the Dirichlet boundary conditionon the upper hemisphere.} \label{Fig:LBUpperHemisphereDir}
\end{figure}

\begin{figure}[!htb]
\centering{
\includegraphics[width=12cm]{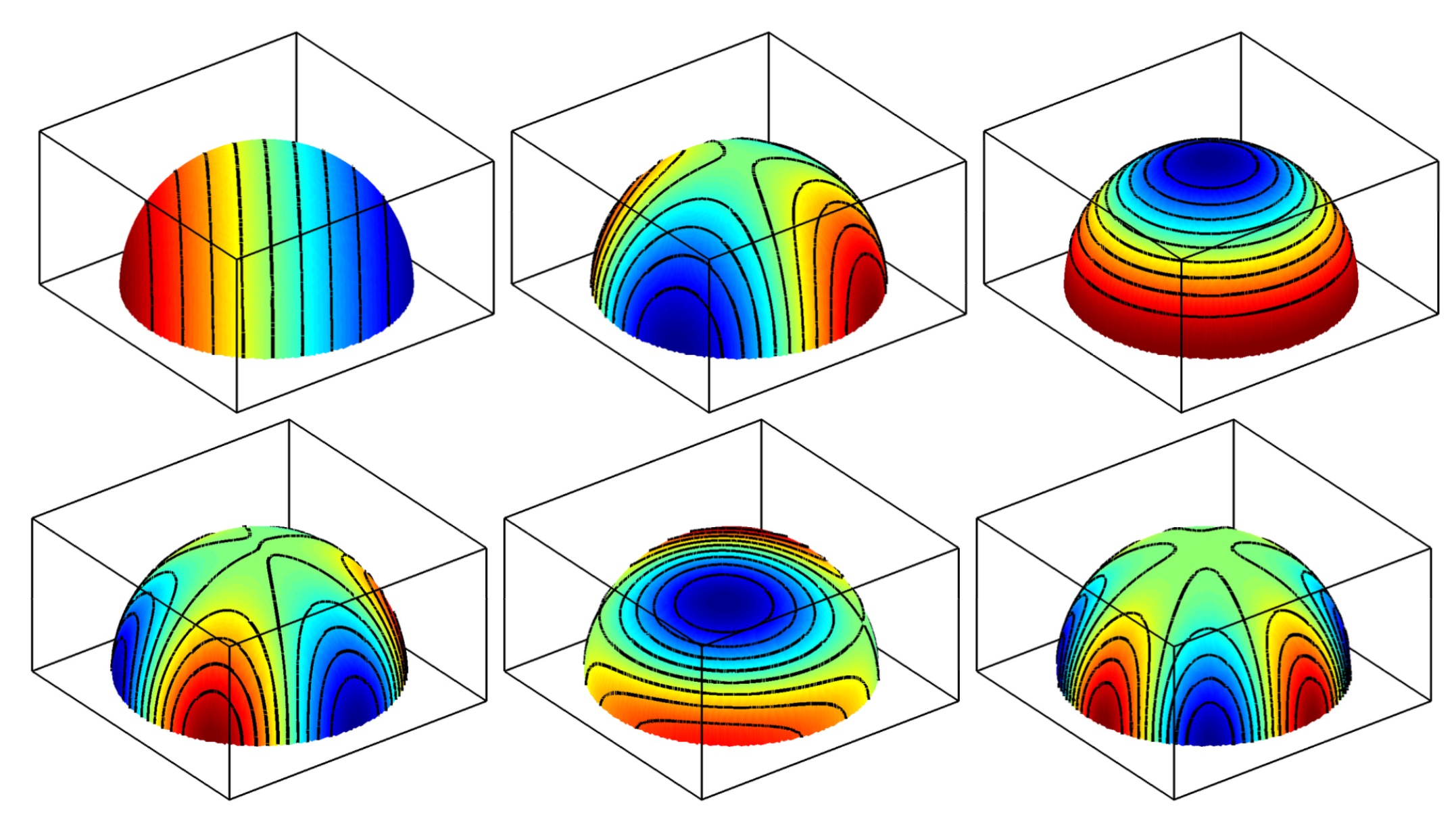}
}
\caption{First few eigenfunctions of the LB operator with the Neumann boundary condition on the upper hemisphere.} \label{Fig:LBUpperHemisphereNeu}
\end{figure}

The proposed MVGD approach can also be applied to open surfaces. We compute the eigenvalues and their corresponding eigenfunctions on an upper hemisphere with the Dirichlet boundary condition $\nu|_{\partial \Sigma}=0$ and the Neumann boundary condition $\partial_{\mathbf{n}} \nu |_{\partial \Sigma}=0$ imposed on the unit circle on the $x$-$y$ plane. Solutions to these problems can be analytically computed. The multiplicities for the eigenvalue $\lambda_n$ in both problems are given by $n$ and their values are given by $\lambda_n=n(n+1)$ and $\lambda_n=n(n-1)$, respectively. Some eigenfunctions are plotted in Fig \ref{Fig:LBUpperHemisphereDir} and Fig \ref{Fig:LBUpperHemisphereNeu}, respectively. Table \ref{table:EigTableUpperSphereDir} shows the $L^{\infty}$ error, $E_n^{\infty}$ as defined above, for $\lambda_5$ and $\lambda_{13}$ for the Dirichlet problem. The $L^{\infty}$ errors in the eigenvalues from the Neumann problem are shown in Table \ref{table:EigTableUpperSphereNeu}. As we increase the number of sampling points, the errors in these numerical solutions are reduced approximately linearly.

\begin{table}
\begin{center}
\begin{tabular}{|c||c|c|c|c|c|}
\hline
Sample size & 500 & 1000 & 2000 & 4000 & 8000 \\
\hline
Boundary points & 20 & 40 & 80 & 160 & 320 \\
\hline
\hline
$\lambda_5=30$ & 1.27 & 0.46 & 0.54 & 0.22 & 0.14 \\
\hline
$\lambda_{13}=182$ & 10.02 & 3.76 & 2.16 & 1.57 & 0.85 \\
\hline
\end{tabular}
\captionof{table}{ $L^{\infty}$ ($\times10^{-2}$) errors in the eigenvalues of the LB operator with the Dirichlet boundary condition on the upper hemisphere.}
\label{table:EigTableUpperSphereDir}
\end{center}
\end{table}

\begin{table}
\begin{center}
\begin{tabular}{|c||c|c|c|c|c|}
\hline
Sample size & 500 & 1000 & 2000 & 4000 & 8000 \\
\hline
Boundary points & 20 & 40 & 80 & 160 & 320 \\
\hline
\hline
$\lambda_5=30$ & 1.71 & 0.86 & 0.38 & 0.18 & 0.12 \\
\hline
$\lambda_{13}=182$ & 7.86 & 4.75 & 2.3 & 1.07 & 0.57 \\
\hline
\end{tabular}
\captionof{table}{$L^{\infty}$ ($\times10^{-2}$) errors in the eigenvalues of the LB operator with the Neumann boundary condition on the upper hemisphere.}
\label{table:EigTableUpperSphereNeu}
\end{center}
\end{table}

\section{Conclusion}
In this paper, we have proposed the Modified Virtual Grid Difference (MVGD) method for discretizing the LB operator on manifolds sampled by point clouds. The discretization is very simple to implement. The discretized system can achieve diagonal dominance and can be efficiently solved by many well developed iterative methods which then leads to a computationally very efficient method for solving the LB equation on manifolds. As a future application, we will extend this proposed LB discretization to solve differential equations on moving interfaces.

\section*{Acknowledgments}
The work of Leung was supported in part by the Hong Kong RGC under grants 16303114 and 16309316. Zhao's research is partially supported by the NSF grant DMS-1418422. The authors would also like to thank J. Liang and R. Lai for their helps on implementing the least squares approach as proposed in \cite{liazha13} and also the local mesh method as proposed in \cite{lailiazha13}.

\bibliographystyle{plain}

\begin{thebibliography}{}

\end{thebibliography}


\begin{thebibliography}{10}

\bibitem{alexa2003computing}
Marc Alexa, Johannes Behr, Daniel Cohen-Or, Shachar Fleishman, David Levin, and
  Claudio~T Silva.
\newblock Computing and rendering point set surfaces.
\newblock {\em IEEE Transactions on visualization and computer graphics},
  9:3--15, 2003.

\bibitem{belniy05}
M.~Belkin and P.~Niyogi.
\newblock Towards a theoretical foundation for {L}aplacian-based manifold
  methods.
\newblock {\em COLT}, pages 486--500, 2005.

\bibitem{belsunwan09}
M.~Belkin, J.~Sun, and Y.~Wang.
\newblock Constructing {L}aplace operator from point clouds in rd.
\newblock In {\em Proceedings of the twentieth Annual ACM-SIAM Symposium on
  Discrete Algorithms}, pages 1031--1040. Society for Industrial and Applied
  Mathematics, 2009.

\bibitem{coilaf06}
R.R. Coifmam and S.~Lafon.
\newblock Diffusion maps.
\newblock {\em Appl. Comput. Harmon. Anal.}, 21(1):5--30, 2006.

\bibitem{dzi88}
G.~Dziuk.
\newblock Finite elements for the beltrami operator on arbitrary surfaces.
\newblock {\em In S. Hilde-brandt and R. Leis, editors, Partial differential
  equations and calculus of variations}, 1357:142--155, 1988.

\bibitem{fleishman2005robust}
Shachar Fleishman, Daniel Cohen-Or, and Cl{\'a}udio~T Silva.
\newblock Robust moving least-squares fitting with sharp features.
\newblock In {\em ACM transactions on graphics (TOG)}, volume~24, pages
  544--552. ACM, 2005.

\bibitem{frenie80}
R.~Franke and G.~Nielson.
\newblock Smooth interpolation of large sets of scattered data.
\newblock {\em Internat. J. Numer. Methods Engrg.}, 15(11):1691--1704, 1080.

\bibitem{gon10}
{\'A}lvaro Gonz{\'a}lez.
\newblock Measurement of areas on a sphere using {F}ibonacci and
  latitude--longitude lattices.
\newblock {\em Mathematical Geosciences}, 42(1):49--64, 2010.

\bibitem{GSC95}
T.N~T Goodman, H.B. Said, and L.H~T Chang.
\newblock Local derivative estimation for scattered data interpolation.
\newblock {\em Applied Mathematics and Computation}, 68:41--50, 1995.

\bibitem{honleuzha14}
S.~Hon, S.~Leung, and H.-K. Zhao.
\newblock A cell based particle method for modeling dynamic interfaces.
\newblock {\em J. Comp. Phys.}, 272:279--306, 2014.

\bibitem{kol98}
R.~Kolluri.
\newblock Provably good moving least squares.
\newblock {\em ACM Transactions on Algorithms}, 4(2):1--25, 2008.

\bibitem{krsek1998algorithms}
Pavel Krsek, Gabor Luk{\'a}cs, and RR~Martin.
\newblock Algorithms for computing curvatures from range data.
\newblock {\em the Mathematics of Surfaces VIII}, pages 1--16, 1998.

\bibitem{lailiazha13}
R.~Lai, J.~Liang, and H.~Zhao.
\newblock A local mesh method for solving pdes on point clouds.
\newblock {\em Inverse Problems and Imaging}, 7(3):737--755, 2013.

\bibitem{lazmon02}
D.~Lazzaro and L.B. Montefusco.
\newblock Radial basis functions for the multivariate interpolation of large
  scattered data sets.
\newblock {\em J. Comput. Appl. Math.}, 140:521--536, 2002.

\bibitem{leulowzha11}
S.~Leung, J.~Lowengrub, and H.K. Zhao.
\newblock A grid based particle method for high order geometrical motions and
  local inextensible flows.
\newblock {\em J. Comput. Phys.}, 230:2540--2561, 2011.

\bibitem{leuzha0802}
S.~Leung and H.K. Zhao.
\newblock A grid-based particle method for evolution of open curves and
  surfaces.
\newblock {\em J. Comput. Phys.}, 228:7706--7728, 2009.

\bibitem{leuzha0801}
S.~Leung and H.K. Zhao.
\newblock A grid based particle method for moving interface problems.
\newblock {\em J. Comput. Phys.}, 228:2993--3024, 2009.

\bibitem{leuzha10}
S.~Leung and H.K. Zhao.
\newblock Gaussian beam summation for diffraction in inhomogeneous media based
  on the grid based particle method.
\newblock {\em Communications in Computational Physics}, 8:758--796, 2010.

\bibitem{lishisun15}
Z.~Li, Z.~Shi, and J.~Sun.
\newblock {Point Integral Method for Solving Poisson-type Equations on
  Manifolds from Point Clouds with Convergence Guarantees}.
\newblock {\em arXiv preprint arXiv:1409.2623}, 2015.

\bibitem{llwz12}
J.~Liang, R.~Lai, T.W. Wong, and H.~Zhao.
\newblock Geometric understanding of point clouds using {L}aplace-{B}eltrami
  operator.
\newblock {\em Computer Vision and Pattern Recongnition (CVPR)}, pages
  214--221, 2012.

\bibitem{liazha13}
J.~Liang and H.~Zhao.
\newblock Solving partial differential equations on point clouds.
\newblock {\em SIAM J. on Scientific Computing}, 35(3):A1461--A1486, 2013.

\bibitem{liuleu13}
J.~Liu and S.~Leung.
\newblock A splitting algorithm for image segmentation on manifolds represented
  by the grid based particle method.
\newblock {\em J. Sci. Comput.}, 56(2):243--266, 2013.

\bibitem{liu2008local}
Ligang Liu, Lei Zhang, Yin Xu, Craig Gotsman, and Steven~J Gortler.
\newblock {\em A local/global approach to mesh parameterization}, volume~27.
\newblock Wiley Online Library, 2008.

\bibitem{magsolriv07}
E.~Magid, O.~Soldea, and E.~Rivlin.
\newblock A comparison of gaussian and mean curvature estimation methods on
  triangular meshes of range image data.
\newblock {\em Comput. Vis. Image Underst.}, 107(3):139--159, 2007.

\bibitem{medvelfig03}
B.~Mederos, L.~Velho, and L.H. De~Figueiredo.
\newblock Moving least squares multiresolution surface approximation.
\newblock {\em Computer Graphics and Image Processing, 2003 SIBGRAPI 2003. XVI
  Brazilian Symposium on 2003}, pages 19--26, 2003.

\bibitem{meek2000surface}
Dereck~S Meek and Desmond~J Walton.
\newblock On surface normal and gaussian curvature approximations given data
  sampled from a smooth surface.
\newblock {\em Computer Aided Geometric Design}, 17:521--543, 2000.

\bibitem{mdsb03}
M.~Meyer, M.~Desbrun, P.~Schroder, and Alan~H Barr.
\newblock Discrete differential-geometry operators for triangulated
  2-manifolds.
\newblock {\em Visualization and Mathematics III}, pages 35--57, 2003.

\bibitem{pauly2003shape}
Mark Pauly, Richard Keiser, Leif~P Kobbelt, and Markus Gross.
\newblock Shape modeling with point-sampled geometry.
\newblock {\em ACM Transactions on Graphics (TOG)}, 22:641--650, 2003.

\bibitem{sinwu15}
A.~Singer and H.-T. Wu.
\newblock Spectral convergence of the connection {L}aplacian from random
  samples.
\newblock {\em Submitted, http://arxiv.org/abs/1306.1587}, 2015.

\bibitem{tau00}
G.~Taubin.
\newblock Geometric signal processing on polygonal meshes.
\newblock {\em EUROGRAPHICS}, 2000.

\bibitem{Willmore93}
T.J. Willmore.
\newblock {\em Riemannian Geometry}.
\newblock New York: Oxford Science Publications, The Clarendon Press, Oxford
  University Press, 1993.

\bibitem{xu04}
G.~Xu.
\newblock Convergent discrete {L}aplace-{B}eltrami operators over triangular
  surfaces.
\newblock {\em Geometric Modeling and Processing}, pages 195--204, 2004.

\end{thebibliography}

\end{document}